\documentclass{amsart}
\usepackage{amsmath,amsthm,bbm}
\usepackage{amsfonts,amssymb}

  \newtheorem{theorem}{Theorem}[section]
  \newtheorem{lemma}[theorem]{Lemma}
  
  \newtheorem{proposition}[theorem]{Proposition}
  \newtheorem{definition}[theorem]{Definition}
  
  \newtheorem{remark}[theorem]{Remark}
  
  \newcounter{figures}[section]

\newcommand{\beq}{\begin{equation}}
\newcommand{\eeq}{\end{equation}}

\newcommand\norm[1]{\|#1\|}
\newcommand\nnorm[1]{\Big\|#1\Big\|}

\newcommand{\R}{\mathbb{R}}
\newcommand{\Rd}{\R^d}

\newcommand{\X}{\mathcal{X}}

  \def\cN{\mathcal{N}}

  \def\cD{\mathcal{D}}
  
  \def\cO{\mathcal{O}}

  \def\ip#1{{\langle #1 \rangle}}

\newcommand{\ONE}{{\mathbbm 1}}
\newcommand{\tONE}{{\tilde{\mathbbm 1}}}

\def\supp{\rm{supp}\, }

\newcommand{\hn}{h_n^{(\alpha,\beta)}}

  \newcommand{\cX}{{\mathcal X}}
  \newcommand{\cM}{{\mathcal M}}

\def\a{{\alpha}}
\def\b{{\beta}}

\newcommand{\wh}{\widehat}
\newcommand{\ha}{\widehat a}
\newcommand{\hb}{\widehat b}

\newcommand\Bes[3]{B^{#1#3}_{#2}}
\newcommand\bes[3]{b^{#1#3}_{#2}}
\newcommand\F[3]{F^{#1#3}_{#2}}
\newcommand\f[3]{f^{#1#3}_{#2}}

\newcommand\BBes[3]{\mathcal{B}^{#1#3}_{#2}}
\newcommand\bbes[3]{\mathbf{b}^{#1#3}_{#2}}
\newcommand\FF[3]{\mathcal{F}^{#1#3}_{#2}}
\newcommand\ff[3]{\mathbf{f}^{#1#3}_{#2}}

\newcommand\BB[2]{\mathcal{B}^{#1}_{#2}}

\newcommand{\eps}{{\varepsilon}}

\newcommand{\w}{w}
\def\W{{\mathcal W}}

\def\cX{{\mathcal X}}
\newcommand{\PP}{{\mathsf {P}}}
\def\a{{\alpha}}
\def\b{{\beta}}
\def\t{{\theta}}

 \def\CO{{\mathcal O}}

 \def\RR{{\mathbb R}}

 \def\supp{\operatorname{supp}}
 
\newcommand{\wt}{\widetilde}

 \def\CO{{\mathcal O}}

 \def\RR{{\mathbb R}}

 \def\supp{\operatorname{supp}}
 
\def\wt{\widetilde}
\def\wh{\widehat}
\def\omab{{\W}}
\def\eps{{\varepsilon}}
\def\dist{{\rm dist\,}}
\def\CJ{{\mathcal J}}
\def\CK{{\mathcal K}}
\def\Ker{{L}}
\def\Lp{{p}}

\def\LL2{{L^2(\w)}}
\def\L2{{2}}
\def\Lpp{{L^p(\w)}}

\def\Linfty{{\infty}}
\def\bC{{\mathbb C}}
\def\Ln{{\Ker_n}}

\def\CC{{\mathbb C}}

\def\a{{\alpha}}
\def\ph{{\varphi}}
\renewcommand{\lesssim}{{\le c}}
\def\cc{{\rm c}}
\def\cY{{\mathcal Y}}

\def\dm{{dm_{\a, \b}}}

\begin{document}

\title{Jacobi decomposition of weighted Triebel-Lizorkin and Besov spaces}

\author{George Kyriazis}
\address{Department of Mathematics and Statistics,
University of Cyprus, 1678 Nicosia, Cyprus}
\email{kyriazis@ucy.ac.cy}

\author{Pencho Petrushev}
\address{Department of Mathematics,
University of South Carolina, Columbia, SC 29208, USA}
\email{pencho@math.sc.edu}


\author{Yuan XU}
\address{Department of Mathematics,
University of Oregon, Eugene, Oregon, 97403, USA}
\email{yuan@math.uoregon.edu}

\date{\today}


\begin{abstract}
The Littlewood-Paley theory is extended to weighted  spaces of
distributions on $[-1,1]$ with Jacobi weights
$
\w(t)=(1-t)^\alpha(1+t)^\beta.
$
Almost exponentially localized polynomial elements (needlets)
$\{\varphi_\xi\}$, $\{\psi_\xi\}$ are constructed
and,
in complete analogy with the classical case on $\RR^n$,
it is shown that weighted Triebel-Lizorkin and Besov spaces
can be characterized by the size of the needlet coefficients
$\{\ip{f,\varphi_\xi}\}$
in respective sequence spaces.
\end{abstract}

\maketitle

\maketitle

\section{Introduction}\label{introduction}
\setcounter{equation}{0}

The $\varphi$-transform of Frazier and Jawerth \cite{FJ1, FJ2, FJW}
is a powerful tool for decomposition of spaces of functions or
distributions on $\RR^n$.
Our goal in this paper is to develop similar tools for decomposition of
weighted spaces of distributions on $[-1,1]$ with Jacobi weights
\begin{equation}\label{wab}
\w(x):=\w_{\a,\b}(x):=(1-x)^\alpha(1+x)^\beta,
\quad \a,\b > -1/2.
\end{equation}
We will build upon the elements constructed in \cite{PX1}
and termed {\em needlets}.
The targeted spaces are weighted Triebel-Lizorkin and Besov spaces
on $[-1,1]$.

The main vehicle in constructing our building blocks will be the
classical Jacobi polynomials
$\{P_n^{(\a,\b)}\}_{n=0}^\infty$,
which form an orthogonal basis for $L^2(\w):=L^2([-1, 1], \w)$ and
are normalized by
$P_n^{(\a,\b)}(1)= \binom{n+\a}{n}$
\cite{Sz}.
In particular,
\begin{equation}\label{h-ab}
\int_{-1}^1 P_n^{(\alpha,\beta)}(x)P_m^{(\alpha,\beta)}(x)\w(x) dx=\delta_{n,m}\hn,
\end{equation}
where $\hn\sim n^{-1}$ with constants of equivalence depending only on $\a$ and $\b$.
Then the normalized Jacobi polynomials $\PP_n(x)=\PP_n^{(\a,\b)}(x)$, defined by
\begin{equation}\label{def-Jacobi}
\PP_n(x):=(\hn)^{-1/2}P_n^{(\alpha,\beta)}(x), \quad n=0, 1, \dots,
\end{equation}
form an orthonormal basis for $L^2(\w)$, where the inner product
is defined by
\begin{equation}\label{def-inn-prod}
\ip{f,g}:=\int_{-1}^1f(x) \overline{g(x)}\w(x)dx.
\end{equation}
Consequently, for every $f\in L^2(\w)$
\beq\label{J-decomposition}
f=\sum_{n=0}^\infty a_n(f) \PP_n
\quad\text{with}\quad
a_n(f):=\ip{f,\PP_n}.
\eeq
Then the kernel of the $n$th partial sum operator is
\beq\label{Bad-kernel}
K_n(x,y):=\sum_{\nu=0}^{n} \PP_\nu(x)\PP_\nu(y).
\eeq

Our construction of needlets relies on the fundamental
fact \cite{PX1} that if the coefficients on the right in (\ref{Bad-kernel})
are ``smoothed out" by sampling a compactly supported $C^\infty$ function,
then the resulting kernel has nearly exponential localization
around the main diagonal $y=x$ in $[-1, 1]^2$.
To be more specific, let
\begin{equation} \label{def-Ln}
L_n(x,y):= \sum_{j=0}^\infty \ha\Big(\frac{j}{n}\Big) \PP_j(x)\PP_j(y)
\end{equation}
with $\ha$ admissible in the sense of the following definition:


\begin{definition} \label{defn:admissible}
A function $\ha \in C^\infty[0, \infty)$
is said to be admissible of type

$(a)$ if
$\supp \ha \subset [0,2]$ and $\ha(t)=1$ on $[0, 1]$,
and of type

$(b)$ if $\supp \ha \subset [1/2,2]$.
\end{definition}
\noindent
As a companion to the weight $\w(x)$ we introduce the quantity
\begin{equation}\label{wn_weight}
\W(n; x)= \W_{\a,\b}(n; x):= (1-x+ n^{-2})^{\a+1/2}(1+x+ n^{-2})^{\b+1/2}.
\end{equation}
We will also need the distance on $[-1,1]$ defined by
\begin{equation} \label{def.dist}
d(x, y):=|\arccos x-\arccos y|.
\end{equation}
Now one of the main results from \cite{PX1} can be stated as follows:
Let $\ha$ be admissible.
Then for any $\sigma>0$ there is a constant $c_\sigma > 0$ depending only on
$\sigma$, $\a$, $\b$, and $\ha$ such that
\begin{equation} \label{Lbound1}
|L_n(x, y)|
\le c_\sigma \frac{n}{\sqrt{\W(n;x)}
\sqrt{\W(n;y)} (1+nd(x, y))^{\sigma}},
\quad x, y\in [-1, 1].
\end{equation}

The kernels $L_n(x, y)$ are the main ingredient in constructing needlet systems here.
Our construction utilizes a semi-discrete Calder\'{o}n type decomposition
combined with discretization using the Gaussian quadrature formula
(see \S\ref{def-needlets}).
Earlier in \cite{NPW} a similar scheme has been used  for the construction of frames
on the sphere.

Denoting by $\{\ph_\xi\}_{\xi\in\cX}$ and $\{\psi_\xi\}_{\xi\in\cX}$
the constructed analysis and synthesis needlet systems, indexed by
a multilevel set $\cX=\cup_{j=0}^\infty \cX_j$,
we show that every distribution $f$ on $[-1, 1]$ ($f\in\cD'$)
has the representation
$$
f=\sum_{\xi \in \cX}\langle f, \ph_\xi \rangle \psi_\xi.
$$

In this article we use the needlets to characterize two scales of
weighted Triebel-Lizorkin (F-space) and Besov spaces (B-spaces) on $[-1, 1]$
defined via Jacobi expansions.
The idea of using orthogonal or spectral decompositions for
introducing Triebel-Lizorkin and Besov spaces is natural and well known,
see \cite{Pee, T1}.
To be more precise, let
$$
\Phi_0(x, y) := \PP_0(x)\PP_0(y)
\quad\mbox{and}\quad
\Phi_j(x, y) := \sum_{\nu=0}^\infty
\ha \Big(\frac{\nu}{2^{j-1}}\Big)\PP_\nu(x)\PP_\nu(y),
\quad j\ge 1,
$$
where $\ha$ is admissible of type (b) (see Definition~\ref{defn:admissible})
and $\ha>0$ on $[3/5, 5/3]$.

The first scale of $F$-spaces $\F spq$ with
$s \in \R$, $0<p<\infty$, $0<q\le\infty$,
is defined (\S\ref{Tri-Liz})
as the space of all distributions $f$ on $[-1, 1]$ such that
$$ 
\norm{f}_{\F spq}:=\nnorm{\biggl(\sum_{j=0}^{\infty}
(2^{sj}|\Phi_j\ast f(\cdot)|)^q\biggr)^{1/q}}_{\Lpp} <\infty.
$$ 
We define a second scale of $F$-spaces $\FF spq$
(\S\ref{Tri-Liz-2}) as the space of all $f\in\cD'$ such that
$$
\|f\|_{\FF spq} :=\nnorm{\biggl(\sum_{j=0}^{\infty}
\Big[2^{sj}\W(2^j;\cdot)^{-s}|\Phi_j\ast f(\cdot)|\Big]^q\biggr)^{1/q}}_{\Lpp}
<\infty.
$$
(For the definition of $\Phi_j\ast f$, see (\ref{convolution1}).)
The corresponding scales of weighted Besov spaces
$\Bes spq$ (see \cite{RS, T1}) and $\BBes spq$
with $s \in \R$, $0<p, q\le\infty$,
are defined
(\S\ref{Besov}-\ref{Besov2}) via the (quasi-)norms
$$
\|f\|_{\Bes spq} :=
\Big(\sum_{j=0}^\infty \Big(2^{s j}
\|\Phi_j\ast f\|_{\Lpp}\Big)^q\Big)^{1/q}
$$
and
$$
\|f\|_{\BBes spq} :=
\Big(\sum_{j=0}^\infty \Big[2^{s j}
\|\W(2^j;\cdot)^{-s}\Phi_j*f(\cdot)\|_{\Lpp}\Big]^q\Big)^{1/q}.
$$
To some extent the second scales of F- and B-spaces are more natural
than the first scales since they scale (embed) correctly with respect
to the smoothness parameter~$s$
(see \S\ref{Tri-Liz-2}, \S\ref{Besov2} for details).
Also, the second scale of B-spaces provides the smoothness spaces
of nonlinear n-term approximation from needlets (\S\ref{Nonlin-app}).

One of our main results (\S\ref{Tri-Liz}) shows that for all indices
the weighted Triebel-Lizorkin spaces $\F spq$
can be characterized in terms of the size of the needlet coefficients,
namely,
$$
\norm{f}_{\F spq}
\sim \nnorm{\biggl(\sum_{j=0}^\infty 2^{sjq}
\sum_{\xi\in \cX_j}|\ip{f, \varphi_\xi}\psi_\xi(\cdot)|^q\biggr)^{1/q}}_{\Lpp}.
$$
The needlet characterization of the Besov spaces $\Bes spq$
(\S\ref{Besov}) takes the form
$$
\norm{f}_{\Bes spq}
\sim \Big(\sum_{j=0}^\infty2^{sjq}
\Bigl[\sum_{\xi\in \cX_j}
\norm{\ip{f,\varphi_\xi}\psi_\xi}_{\Lp}^p\Bigr]^{q/p}\Bigr)^{1/q}.
$$
Characterizations of similar nature are obtained for the second scales of
weighted Triebel-Lizorkin and Besov spaces $\FF spq$ and $\BBes spq$
(see \S\ref{Tri-Liz-2}, \S\ref{Besov2}).
Using $\Lpp$ multipliers we show that the space $\F 0p2 = \FF 0p2$
can be identified as $\Lpp$ for $1<p<\infty$.

Atomic and molecular decomposition of weighted
Triebel-Lizorkin and Besov spaces can be developed using the approach
of Frazier and Jawerth \cite{FJ1, FJ2}.
This enables one to make the connection between the weighted F-spaces
and the weighted Hardy spaces on $[-1, 1]$ (see \cite{CW, MS}).
To prevent this paper from exceeding some reasonable size
we leave the atomic and molecular decompositions for elsewhere.

It is an open problem to extend the results from this article to dimensions
$d>1$. The missing key element is the nearly exponential localization of
kernels of type (\ref{def-Ln}) in the multivariate case.

The rest of the paper is organized as follows.
In \S\ref{preliminaries}, some auxiliary facts are given, including
localized and reproducing polynomial kernels,
Gaussian quadrature, the maximal inequality,
and basics of distributions on $[-1, 1]$.
In \S\ref{def-needlets}, we construct the needlets and show some of their
properties.
The first and second scales of weighted Triebel-Lizorkin spaces are defined
and characterized via neadlets in \S\ref{Tri-Liz} and \S\ref{Tri-Liz-2},
respectively, while
the first and second scales of Besov spaces are defined and characterized
via needlets in \S\ref{Besov} and \S\ref{Besov2}.
In \S\ref{Nonlin-app}, Besov spaces are applied to weighted nonlinear
approximation from needlets; a Jackson theorem is proved.
Section \ref{appendix} is an appendix, where the proofs of some
statements are given.

Throughout the paper we use the following notation:
$$
\|f\|_\Lp:=\Big(\int_{-1}^1|f(x)|^p\w(x)dx\Big)^{1/p},
\quad 0<p<\infty,
\quad\mbox{and}
\quad
\|f\|_\Linfty:=\sup_{x\in[-1,1]}|f(x)|.
$$
For a measurable set $E\subset [-1, 1]$, we set $\mu(E):=\int_E\w(y)\, dy$;
$\ONE_E$ is the characteristic function of $E$
and $\tONE_E:=|\mu(E)|^{-1/2}\ONE_E$ is the $L^2(\w)$ normalized characteristic
function of $E$.
Also, $\Pi_n$ denotes the set of all univariate algebraic polynomials
of degree $\le n$.
Positive constants are denoted by $c$, $c_1, c_*, \dots$  and they
may vary at every occurrence. The notation $A\sim B$
means $c_1A\le B\le c_2 A$.

\section{Preliminaries}\label{preliminaries}
\setcounter{equation}{0}

\subsection{Localized kernels induced by Jacobi polynomials}\label{Local-kernels}

To a large extent our development in this paper relies on the nearly exponential
localization (\ref{Lbound1}) of kernels $L_n(x, y)$ of the form (\ref{def-Ln})
with admissible $\ha$, established in \cite{PX1}.
To~avoid some potential confusion, we note that
the inner product in \cite{PX1} is defined by $
\ip{f,g}:=c_{\a,\b}\int_{-1}^1f(x) \overline{g(x)}\w(x)dx $ with
$c_{\a,\b}^{-1}:=\int_{-1}^1\w(x)dx$ and as a
result $L_n(x, y)$ from (\ref{def-Ln})
is a constant multiple of $L^{\a,\b}(x,y)$ from \cite{PX1}.
A similar remark applies to the constants $h_n^{(\a,\b)}$
from (\ref{h-ab}) and \cite{PX1}.

The proof of estimate (\ref{Lbound1}) (see \cite{PX1}) is
based on the almost exponential localization of the univariate polynomial:
\begin{equation}\label{def.Ln}
L_n^{\a,\b}(x) := \sum_{j=0}^\infty \wh a\left(\frac{j}{n}\right)
  \left(h_j^{(\a,\b)}\right)^{-1} P_j^{(\a,\b)}(1)P_j^{(\a,\b)}(x).
\end{equation}


\begin{theorem}\cite{BD, PX1}\label{thm:localL-x}
Assume that $\a \ge \b > -1/2$ and let $\ha$ be admissible.
Then for every $k\ge 1$ there exists a constant $c_k>0$ depending only
on $k$, $\a$, $\b$, and $\wh a$ such that
\begin{equation}\label{est.L-n}
|L_n^{\a,\b}(\cos \theta)| \le c_k \frac{n^{2\a+2}}{(1+n \theta)^{k+\a-\b}},
  \quad 0 \le \theta \le \pi.
\end{equation}
The dependence of $c_k$ on $\ha$ is of the form
$c_k=c(\a, \b, k)\max_{1\le\nu\le k}\|\ha^{(\nu)}\|_{L^1}$.
\end{theorem}
This estimate was proved in
\cite{PX1} with $\ha$ admissible of type (b) and
in \cite{BD} with $\ha$ admissible of type (a)
(for a proof, see also~\cite{PX2}).

\smallskip

In \cite[Proposition 1]{PX1} it is shown that (\ref{Lbound1}) yields
the following upper bound for the weighted $L^p$ integrals of $|L_n(x,y)|$:
\begin{equation} \label{est-Lp-int}
\int_{-1}^1 |L_n (x,y)|^p \w(y) dy \le
c\Big(\frac{n}{\W(n;x)}\Big)^{p-1},  \quad -1\le x\le 1, \quad 0 < p<\infty.
\end{equation}

The next theorem shows that in a sense the kernel $\Ker_n(x, y)$
from (\ref{def-Ln}) 
is Lip 1 in $x$ (and $y$).


\begin{theorem}\label{thm:Lip}
Let $\a,\b > -1/2$. Suppose $\ha$ is admissible
and $\sigma>0$ is an arbitrary constant.
If $x, y, z, \xi \in [-1, 1]$,
$d(x, \xi) \le c_*n^{-1}$ and $d(z, \xi) \le c_*n^{-1}$
with $n\ge 1$, $c_*>0$, then
\begin{equation}\label{Lip}
|\Ln(x,y) - \Ln(\xi,y)| \le c_\sigma \frac{n^2 d(x, \xi)}
 {\sqrt{\W(n; y)}\sqrt{\W(n; z)} (1+nd(y, z))^{\sigma}},
\end{equation}
where $c_\sigma>0$ depends only on
$\sigma$, $\a$, $\b$, $c_*$, and $\ha$.
\end{theorem}
The proof of this theorem is given in the appendix.

Lower bound estimates for the integrals of $|L_n(x,y)|$ are nontrivial
and will be vital for our further development.

\begin{proposition}
Let $\ha$ be admissible and
$|\ha(t)| \ge c > 0$ for $t \in [3/5,5/3]$. Then
\begin{equation} \label{est-Lp-norm}
\int_{-1}^1 |L_n (x,y)|^2 \w(y) dy \ge c
n\W(n;x)^{-1},  \quad -1\le x\le 1.
\end{equation}
\end{proposition}

\noindent
{\bf Proof.}
By the definition of $L_n(x, y)$ in (\ref{def-Ln})
and the orthogonality of the Jacobi polynomials, it follows that
$$
\int_{-1}^1 |L_n(x,y)|^2 \w(y) dy \ge c
\sum_{k = [n/2] }^{2 n} |\ha(k/n)|^2 [\PP_k^{(\a,\b)}(x)]^2.
$$
Since $|\wh a(t)| \ge c > 0$ for $t \in [3/5,5/3]$ and
$\PP_k(x) = (h_k^{(\a,\b)})^{-1/2} P_k^{(\a,\b)}(x)
\sim k^{1/2} P_k^{(\a,\b)}(x)$,
it suffices to prove that
$$
\sum_{ k = [3 n/5] }^{[5 n/3]} [P_k^{(\a,\b)}(x)]^2
\ge c\W(n;x)^{-1},
\quad c>0,
$$
which is established in the following proposition.
$\qed$


\begin{proposition}\label{prop:lower-bound}
If $\alpha, \beta > -1$ and $\eps>0$, then
\begin{equation} \label{lowerbd}
\Lambda_n(x):= \sum_{ k = n }^{n+[\eps n]} [P_k^{(\a,\b)}(x)]^2
\ge c\W(n;x)^{-1},
\quad x\in [-1, 1],\quad n\ge 1/\eps,
\end{equation}
where $c > 0$ depending only on $\alpha, \beta$, and $\eps$.
\end{proposition}

This proposition is nontrivial and its proof is given in the appendix.


\subsection{Reproducing kernels and best polynomial approximation}
\label{kernels}
We let\\ $E_n(f)_p$ denote the best approximation of $f \in \Lpp$ from
$\Pi_n$, i.e.
\begin{equation}\label{def-En}
E_n(f)_p := \inf_{g \in \Pi_n}\|f-g\|_\Lp.
\end{equation}

To simplify our notation we introduce the following ``convolution":
For functions $\Phi: [-1, 1]^2 \to \bC$ and $f: [-1, 1] \to \bC$,
we write
\begin{equation}\label{convolution}
\Phi*f(x) := \int_{-1}^1 \Phi(x, y)f(y)\w(y)\,dy.
\end{equation}


\begin{lemma}\label{lem:Ker-n}
Suppose $\ha$ is admissible of type $(a)$
and let $L_n(x, y)$ be the kernel defined in $(\ref{def-Ln})$.

$(i)$
$\Ker_n(x, y)$ is a symmetric reproducing kernel for $\Pi_n$,
i.e.
$\Ker_n*g =g$ for $g \in \Pi_n$.

$(ii)$
For any $f \in \Lpp$, $1\le p \le \infty$, we have
$\Ker_n*f \in \Pi_{2n}$,
\begin{equation}\label{Ker-n3}
\|\Ker_n*f\|_\Lp\le c \|f\|_\Lp,
~~~\mbox{and}~~~
\|f-\Ker_n*f\|_\Lp \le cE_n(f)_p.
\end{equation}
\end{lemma}

\noindent
{\bf Proof.}
Part (i) is immediate since
$\ha(\frac{\nu}{n})=1$ for $0 \le \nu \le n$.
The left-hand side estimate in (\ref{Ker-n3}) follows from (\ref{est-Lp-int})
when $p=1$ and $p=\infty$; the general case follows by interpolation.
The right-hand side estimate in (\ref{Ker-n3}) follows by
the left-hand side estimate and (i).
$\qed$

Lemma~\ref{lem:Ker-n} (i) and (\ref{est-Lp-int}) are instrumental in proving
Nikolski type inequalities.


\begin{proposition}\label{Nikolski}
For $0 < q \le p \le \infty$ and $g \in \Pi_n$,
\begin{equation}\label{norm-relation}
\|g\|_p \le cn^{(2+ 2 \min\{0, \max\{\a, \b\} \})(1/q-1/p)}\|g\|_q,
\end{equation}
furthermore, for any $s\in\RR$,
\begin{equation}\label{norm-relation2}
\|\W(n;\cdot)^{s} g(\cdot)\|_p
\le cn^{1/q-1/p}\|\W(n;\cdot)^{s+1/p-1/q} g(\cdot)\|_q.
\end{equation}
\end{proposition}

The proof of this proposition is given in the appendix.


\subsection{Quadrature formula and subdivision of \boldmath $[-1, 1]$}
\label{cubature}

For the construction of our building blocks (needlets)
we will utilize an appropriate Gaussian quadrature formula.
Let $\xi_{j,\nu}=:\cos \theta_{\nu}$, $\nu=1, 2,\dots, 2^{j+1}$, be
the zeros of the Jacobi polynomials $P_{2^{j+1}}^{(\a,\b)}$
ordered so that
$$
0< \theta_{1} < \dots < \theta_{2^{j+1}} < \pi.
$$
It is well known that uniformly
(see \cite{FW} and also (\ref{asymp-zeros})-(\ref{zero}) below)
\beq\label{distr-zeros}
\theta_{\nu+1} - \theta_\nu \sim 2^{-j}
\quad\mbox{and hence} \quad \theta_\nu \sim\nu 2^{-j}.
\eeq
Define now
\begin{equation}\label{def-Xj}
\cX_j:=\{\xi_{j,\nu}: \nu=1, 2,\dots, 2^{j+1}\}, \quad j\ge 0,
\quad \mbox{and set}\quad
\cX:=\cup_{j=0}^\infty\cX_j.
\end{equation}
As is well known \cite{Sz} the zeros of the Jacobi polynomial
$P_{2^{j+1}}^{(\a,\b)}$ serve as knots of the Gaussian quadrature
\begin{equation}\label{quadrature1}
\int_{-1}^1 f(x)\w(x)dx \sim \sum_{\xi \in \cX_j} \cc_\xi f(\xi),
\end{equation}
which is exact for all polynomials of degree at most $2^{j+2}-1$.
Furthermore, the coefficients $\cc_\xi$ are positive and have the asymptotic
\begin{equation}\label{coeff}
\cc_\xi \sim \lambda_{2^{j+1}}(\xi)\sim 2^{-j} \w(\xi)(1-\xi^2)^{1/2}
\sim 2^{-j} \W(2^j; \xi),
\end{equation}
where $\lambda_{2^{j+1}}(t)$ is the Christoffel function and
the constants of equivalence depend only on $\a,\b$
(cf. e.g. \cite{N}).

We next introduce the $j$th level weighted dyadic intervals.
Set as above $\xi_{j,\nu}=:\cos \theta_{\nu}$ and define
\begin{equation}\label{def:Ixi}
I_{\xi_{j, \nu}}:= [(\xi_{j,\nu+1}+\xi_{j,\nu})/2, (\xi_{j,\nu-1}+\xi_{j,\nu})/2],
\quad \nu=2, 3, \dots, 2^{j+1}-1,
\end{equation}
and
\begin{equation}\label{def:Ixi1}
I_{\xi_{j,1}}:= [(\xi_{j,2}+\xi_{j,1})/2, 1],
\quad  I_{\xi_{j,2^{j+1}}}:= [-1, (\xi_{j,2^{j+1}}+\xi_{j,2^{j+1}-1})/2].
\end{equation}
For $\xi\in \cX_j$ we will briefly denote
$I_\xi:=I_{\xi_{j, \nu}}$ if $\xi=\xi_{j, \nu}$.

It follows by (\ref{distr-zeros}) that there exist constants $c_1, c_2 >0$
such that
\begin{equation}\label{I-B}
B_\xi(c_12^{-j}) \subset I_\xi \subset B_\xi(c_22^{-j}),
\end{equation}
where
$B_y(r):= \{x\in [-1, 1]: d(x, y)\le r\}$
with $d(\cdot, \cdot)$ being the distance from (\ref{def.dist}).
Also, it is straightforward to show that
\begin{equation}\label{muI}
\mu(I_\xi):= \int_{I_\xi}\w(x)\, dx \sim 2^{-j}\W(2^j; \xi)\sim \cc_\xi,
\quad \xi\in \cX_j,\quad j\ge 0.
\end{equation}


It will be useful to note that
\begin{align}
\W(n;\cos \theta)&\sim (\sin\theta + n^{-1})^{2 \a +1},
\quad   0\le \theta\le 2\pi/3,\label{wn-sim1}\\
\W(n;\cos \theta)&\sim (\sin\theta + n^{-1})^{2 \b +1},
\quad  \pi/3\le \theta \le \pi.\label{wn-sim2}
\end{align}

The following simple inequality will be instrumental in various proofs
\begin{equation}\label{omega<omega}
\omab(n; x)\le c\omab(n; y)(1+nd(x, y))^{2\max\{\a, \b\}+1},
\quad x, y\in [-1, 1], \quad n\ge 1.
\end{equation}
For the proof see the appendix.

\subsection{The maximal inequality}\label{max}

For every $0<t<\infty$ and $x\in[-1,1]$, we define
\begin{equation}\label{def.max-fun}
\cM_tf(x):=\sup_{I\ni x}\left(\frac1{\mu(I)}\int_I|f(y)|^t\w(y)\,
dy\right)^{1/t},
\end{equation}
where the sup is over all intervals $I\subset[-1,1]$ containing $x$.
It is not hard to see that
$\mu$ is a doubling measure on $[-1,1]$ and hence the general theory of maximal
inequalities applies.
In particular the Fefferman-Stein vector-valued maximal inequality holds
(see \cite{Stein}):
If $0<p<\infty, 0<q\le\infty$ and $0<t<\min\{p,q\}$
then for any sequence of functions $\{f_\nu\}_{\nu=1}^\infty$ on $[-1,1]$,
\begin{equation}\label{max_ineq}
\nnorm{\Bigl(\sum_{\nu=1}^\infty|\cM_tf_\nu(\cdot)|^q\Bigr)^{1/q}}_{\Lp}
\le c\nnorm{\Bigl(\sum_{\nu=1}^\infty| f_\nu(\cdot)|^q\Bigr)^{1/q}}_{\Lp}.
\end{equation}

We need to estimate $(\cM_t\ONE_{I_\xi})(x)$ for the intervals $I_\xi$ from
(\ref{def:Ixi})-(\ref{def:Ixi1}) and other intervals.


\begin{lemma}\label{lem:J-maximal}
Let $\eta \in [0, 1]$ and $0<\eps\le \pi$.
Then for $x\in [-1, 1]$
\begin{equation}\label{J-max1}
(\cM_t \ONE_{B_\eta(\eps)})(x)
\sim \Big(1+\frac{d(\eta, x)}{\eps}\Big)^{-1/t}
     \Big(1+\frac{d(\eta, x)}{\eps+d(\eta, 1)}\Big)^{-(2\alpha+1)/t}
\end{equation}
and hence
\begin{equation}\label{J-max3}
c'\Big(1+\frac{d(\eta, x)}{\eps}\Big)^{-(2\alpha+2)/t}
\le (\cM_t \ONE_{B_\eta(\eps)})(x)
\le c\Big(1+\frac{d(\eta, x)}{\eps}\Big)^{-1/t}.
\end{equation}
Here the constants depend only on
$\a$, $\b$, and $t$.
\end{lemma}

A similar lemma holds for $\eta\in [-1, 0)$.
We relegate the proof of this lemma to the appendix.


\subsection{Distributions on \boldmath $[-1, 1]$}\label{distributions}

Here we give some basic and well known facts about distributions on $[-1, 1]$.
We will use as test functions the set
$\cD:=C^{\infty}[-1,1]$
of all infinitely differentiable complex valued functions on $[-1,1]$,
where the topology is induced by the semi-norms
\beq \label{semi-norms}
\norm{\phi}_k:=\sup_{-1\le t\le 1}|\phi^{(k)}(t)|,
\quad k=0, 1, \dots.
\eeq
Note that the Jacobi polynomials $\{\PP_n\}$ belong to $\cD$.
More importantly, the space $\cD$ of test functions $\phi$
can be completely characterized by the coefficients of their
Jacobi expansions:
$ a_n(\phi):=\langle \phi, \PP_n\rangle
:=\int_{-1}^1\phi(x)\PP_n(x)\w(x)dx.$
Denote
\beq\label{D-norms}
\cN_k(\phi):=\sup_{n\ge 0}\, (n+1)^k |a_n(\phi)|.
\eeq


\begin{lemma}\label{lem:char-D}
$(a)$
$\phi\in\cD$ if and only if
$a_n(\phi) =\cO(n^{-k})$ for all $k$.

$(b)$
For every $\phi\in \cD$ we have
$
\phi=\sum_{n=0}^\infty a_n(\phi)\PP_n,
$
where the convergence is in the topology of $\cD$.

$(c)$
The topology in $\cD$ can be equivalently defined by the norms
$\cN_k(\cdot)$, $k=0, 1, \dots$.
\end{lemma}

\noindent
{\bf Proof.}
If $\phi\in\cD$, then due to the orthogonality of
$\PP_n$ to $\Pi_{n-1}$, we have for $n=1, 2, \dots$
$$
|a_n(\phi)|=|\langle \phi, \PP_n\rangle|=
|\langle \phi-Q_{n-1}, \PP_n\rangle|\le
E_{n-1}(\phi)_2 \le c_kn^{-k}\|\phi^{(k)}\|_{\Linfty},
$$
where $Q_{n-1}\in \Pi_{n-1}$ is the polynomial of best $\LL2$
approximation to $\phi$.
Here we used a simple Jackson estimate for approximation
from algebraic polynomials
($E_n(\phi)_\infty \le c_kn^{-k}\|\phi^{(k)}\|_{\Linfty}$).
Therefore,
$a_n(\phi)=\cO(n^{-k})$ and
$\cN_k(\phi) \le c_k\|\phi\|_k$ for $k=0, 1, \dots$.

On the other hand, by Markov's inequality it follows that
$$
\|\PP_n^{(k)}\|_{L^\infty[-1, 1]}
\le n^{2k}\|\PP_n\|_{L^\infty[-1, 1]}
\le cn^{2k}h_n^{-1/2}P_n(1)
\le cn^{2k+\a+1/2}.
$$
Hence, if $a_n(\phi)=\cO(n^{-k})$ for all $k$,
then
$
\phi^{(k)}=\sum_{n=0}^\infty a_n(\phi)\PP_n^{(k)}
$
with the series converging uniformly and
$$
\|\phi\|_k \le c\sum_{n=0}^\infty|a_n(\phi)|(n+1)^{2k+\a+1/2}
\le c\cN_{2k+[\a+1/2]+1}(\phi), \quad k=0, 1, \dots,
$$
which completes the proof of the lemma.
$\qed$

\medskip

The space $\cD':=\cD'[-1,1]$ of distributions on $[-1,1]$ is defined
as the set of all continuous linear functionals on $\cD$.
The pairing of $f\in \cD'$ and $\phi\in\cD$ will be denoted by
$\langle f, \phi \rangle:= f(\overline{\phi})$, which will be shown
to be consistent with the inner product
$
\langle f, g \rangle :=\int_{-1}^1 f(x)\overline{g(x)}\w(x)dx
$
in $\LL2$.
We will need the representation of distributions from $\cD'$
in terms of Jacobi polynomials.


\begin{lemma}\label{lem:dec-D1}
$(a)$
A linear functional $f$ on $\cD$ is a distribution $(f\in\cD')$
if and only if there exists $k\ge 0$ such that
\begin{equation}\label{D1}
|\langle f, \phi\rangle|\le c_k\cN_k(\phi)
\quad \mbox{for all } \quad \phi \in \cD,
\end{equation}
For $f\in\cD'$, denote $a_n(f):=\langle f, \PP_n\rangle$.
Then for some $k\ge 0$
\begin{equation}\label{D2}
|\langle f, \PP_n\rangle| \le c_k(n+1)^k,
\quad n=0, 1, \dots.
\end{equation}

$(b)$
Every $f\in\cD'$ has the representation
$
f=\sum_{n=0}^\infty a_n(f)\PP_n
$
in distributional sense, i.e.
\begin{equation}\label{D3}
\langle f, \phi\rangle
=\sum_{n=0}^\infty a_n(f)\langle \PP_n, \phi\rangle
=\sum_{n=0}^\infty a_n(f)a_n(\overline{\phi})
\quad\mbox{for all}
\quad \phi\in\cD,
\end{equation}
where the series converges absolutely.
\end{lemma}

\noindent
{\bf Proof.} (a) Part (a) follows immediately
by the fact that the topology in $\cD$ can be defined by the norms $\cN_k(\cdot)$
defined in (\ref{D-norms}).

(b) Using Lemma~\ref{lem:char-D} (b) we get for $\phi\in\cD$,
$$
\langle f, \phi\rangle
=\lim_{N\to\infty}\Big\langle f, \sum_{n=0}^N a_n(\phi)\PP_n\Big\rangle
=\lim_{N\to\infty}\sum_{n=0}^N a_n(\overline{\phi})\langle f, \PP_n\rangle
=\sum_{n=0}^\infty a_n(f)a_n(\overline{\phi}),
$$
where for the last equality we used (\ref{D2}) and the fact that
$a_n(\overline{\phi})$ are rapidly decaying.
$\qed$

\medskip

It is convenient to us to extend the ``convolution"
from (\ref{convolution}) to the case of distributions.


\begin{definition}
Assuming that $f \in \cD'$ and $\Phi: [-1, 1]^2\to\CC$ is such that
$\Phi(x, y)$ belongs to $\cD$ as a function of $y$ $(\Phi(x, \cdot)\in \cD)$,
we define $\Phi*f$ by
\begin{equation}\label{convolution1}
\Phi*f (x) := \langle f, \overline{\Phi(x, \cdot)} \rangle,
\end{equation}
where on the right $f$ acts on $\overline{\Phi(x, y)}$ as a function of $y$.
\end{definition}


\section{Construction of building blocks (Needlets)}\label{def-needlets}
\setcounter{equation}{0}

Following the ideas from \cite{PX1} we next construct two
sequences of companion ``analysis" and ``synthesis" needlets.
Our construction is based on a Calderon type reproducing formula.
Let $\ha,\hb$ satisfy the conditions
\begin{align}
&\quad \ha,\hb\in C^\infty[0, \infty),
\quad \supp \ha,\hb \subset [1/2, 2], \label{ha-hb1}\\
&\quad  |\ha(t)|, |\hb(t)| >c>0, \quad \text{if  } t \in [3/5, 5/3],\label{ha-hb2}\\
&\quad \overline{\ha(t)}\, \hb(t) + \overline{\ha(2t)}\,\hb(2t)  =1,\label{ha-hb3}
\quad\text{if } t \in [1/2, 1].
\end{align}
Hence,
\beq\label{part-unity}
\sum_{\nu=0}^\infty \overline{\ha(2^{-\nu}t )}\, \hb(2^{-\nu}t)=1,
\quad t\in [1,\infty).
\eeq
It is easy to show that
if $\ha$ satisfies (\ref{ha-hb1})-(\ref{ha-hb2}), then
there exists $\hb$ satisfying (\ref{ha-hb1})-(\ref{ha-hb2})
such that (\ref{ha-hb3}) holds true (see e.g. \cite{FJ2}).

\smallskip

Assuming that $\ha$, $\hb$ satisfy (\ref{ha-hb1})-(\ref{ha-hb3}), we define
$\Phi_0(x, y)= \Psi_0(x, y):=\PP_0(x)\PP_0(y)$,
\begin{equation}\label{def.Phi-j}
\Phi_j(x, y) := \sum_{\nu=0}^\infty
\ha\Big(\frac{\nu}{2^{j-1}}\Big)\PP_\nu(x)\PP_\nu(y),
\quad j\ge 1,
\end{equation}
\begin{equation}\label{def-Psi-j}
\Psi_j(x, y) := \sum_{\nu=0}^\infty
\hb\Big(\frac{\nu}{2^{j-1}}\Big)\PP_\nu(x)\PP_\nu(y),
\quad j\ge 1.
\end{equation}
Let $\cX_j$ be the set of knots of quadrature formula \eqref{quadrature1},
defined in (\ref{def-Xj}),
and let $\cc_\xi$ be the coefficients of the same quadrature.
We define the $j$th level {\em needlets} by
\begin{equation}\label{def-needlets1}
\ph_\xi(x) := \cc_\xi^{1/2}\Phi_j(x, \xi)
\quad\mbox{and}\quad
\psi_\xi(x) := \cc_\xi^{1/2}\Psi_j(x, \xi),
\qquad \xi \in \cX_j.
\end{equation}
As in \eqref{def-Xj} we write $\cX := \cup_{j = 0}^\infty \cX_j$,
where equal points from different levels $\cX_j$ are considered
as distinct elements of $\cX$, so that $\cX$ can be used as an index set.
We define the {\em analysis} and {\em synthesis} needlet systems
$\Phi$ and $\Psi$ by
\begin{equation}\label{def-needlets2}
\Phi:=\{\ph_\xi\}_{\xi\in\cX}, \quad \Psi:=\{\psi_\xi\}_{\xi\in\cX}.
\end{equation}

By estimate (\ref{Lbound1}) it follows that the needlets have
nearly exponential localization, namely, for $x\in [-1, 1]$,
\begin{equation}\label{local-Needlets2}
|\Phi_j(\xi, x)|, |\Psi_j(\xi, x)|
\le \frac{c_\sigma 2^{j}}
{\sqrt{\W(2^j; \xi)}\sqrt{\W(2^j; x)}(1+2^{j}d(\xi, x))^\sigma}
\quad \forall \sigma,
\end{equation}
and hence
\begin{equation}\label{local-Needlets21}
|\ph_\xi(x)|, |\psi_\xi(x)|
\le \frac{c_\sigma2^{j/2}}{\sqrt{\W(2^j; x)}(1+2^{j}d(\xi, x))^\sigma}
\quad \forall \sigma.
\end{equation}
Note that $x$ in the term $\sqrt{\omab(2^j; x)}$ above can be replaced by $\xi$
(replacing $c_\sigma$ by a larger constant), namely,
\begin{equation}\label{local-needlets3}
|\ph_\xi(x)|, |\psi_\xi(x)|
\le \frac{c_\sigma2^{j/2}}{\sqrt{\W(2^j; \xi)}(1+2^{j}d(\xi, x))^\sigma}
\quad \forall \sigma.
\end{equation}
This estimate follows by (\ref{local-Needlets21}) and (\ref{omega<omega}).


We will need to estimate the norms of the needlets.
We have for $0<p\le\infty$,
\begin{equation}\label{norm-Needlets}
\|\ph_\xi\|_\Lp \sim \|\psi_\xi\|_\Lp \sim \|\tONE_{I_\xi}\|_\Lp
\sim \Big(\frac{2^j}{\W(2^j; \xi)}\Big)^{1/2-1/p},
\quad \xi\in\cX_j.
\end{equation}
Moreover, there exist constants $c^*, c^\diamond >0$ such that
\begin{equation}\label{norm-Needlets2}
\|\ph_\xi\|_{L^\infty(B_\xi(c^*2^{-j}))},\;
\|\psi_\xi\|_{L^\infty(B_\xi(c^*2^{-j}))}
\ge c^\diamond \Big(\frac{2^j}{\W(2^j; \xi)}\Big)^{1/2},
\end{equation}
where
$
B_\xi(c^*2^{-j}):=\{x\in [-1, 1]: d(\xi, x) \le c^*2^{-j}\}
$
which is an interval.
Notice that if $\ha$, $\hb$ in (\ref{ha-hb1})-(\ref{ha-hb3})
are real valued then by Proposition~\ref{prop:lower-bound}
\begin{equation}\label{Needlet-low}
\ph_\xi(\xi), \psi_\xi(\xi)
\ge c\Big(\frac{2^j}{\W(2^j; \xi)}\Big)^{1/2},
\quad c>0.
\end{equation}
For the proofs of (\ref{norm-Needlets})-(\ref{norm-Needlets2}),
see the appendix.

Our next goal is to establish needlet decompositions
of $\cD'$ and $\Lpp$.


\begin{proposition}\label{prop:needlet-rep}
$(a)$ For $f \in \cD'$, we have
\begin{equation}\label{Needle-rep}
f = \sum_{j=0}^\infty
\Psi_j*\overline{\Phi}_j*f
\quad\mbox{in} ~ \cD'
\end{equation}
and
\begin{equation}\label{needlet-rep1}
f = \sum_{\xi \in \cX}
\langle f, \ph_\xi\rangle \psi_\xi
\quad\mbox{in} ~ \cD'.
\end{equation}

$(b)$ If $f \in \Lpp$, $1\le p \le \infty$, then
$(\ref{Needle-rep})-(\ref{needlet-rep1})$ hold in $\Lpp$.
Moreover, if $1 < p < \infty$, then the convergence in
$(\ref{Needle-rep})-(\ref{needlet-rep1})$ is unconditional.
\end{proposition}

\noindent
{\bf Proof.}
(a) Let $f\in \cD'$.
By (\ref{convolution1}) and Lemma~\ref{lem:dec-D1},
we have
\begin{equation}\label{conv}
\overline{\Phi}_j\ast f
=\sum_{\nu=0}^{2^j} \overline{\ha\Big(\frac{\nu}{2^{j-1}}\Big)}a_\nu(f)\PP_\nu
\end{equation}
and further
\begin{equation}\label{conv-conv}
\Psi_j\ast\overline{\Phi}_j\ast f
=\sum_{\nu=0}^{2^j} \overline{\ha\Big(\frac{\nu}{2^{j-1}}\Big)}
\hb\Big(\frac{\nu}{2^{j-1}}\Big)a_\nu(f)\PP_\nu.
\end{equation}
Now (\ref{Needle-rep}) follows by (\ref{part-unity}) and
Lemmas~\ref{lem:char-D}--\ref{lem:dec-D1}.

Note that
$\Psi_j(x, y)$ and $\Phi_j(x, y)$
are symmetric functions (e.g. $\Psi_j(y, x)=\Psi_j(x, y)$)
and hence
$\Psi_j*\overline{\Phi}_j(x, y)$ is well defined.
Also,
$\Psi_j*(\overline{\Phi}_j*f)=(\Psi_j*\overline{\Phi}_j)*f$.
We observe that
$\Psi_j(x, u)\overline{\Phi_j(y, u)}$
belongs to $\Pi_{2^{j+1}-1}$ as a function of $u$
and apply the quadrature formula from (\ref{quadrature1})
to obtain
\begin{eqnarray*}
\Psi_j*\overline{\Phi}_j(x, y)
&=& \int_{-1}^1 \Psi_j(x, u)\overline{\Phi_j(y, u)}\w(u)du\\
&=&\sum_{\xi\in \cX_j}
\cc_\xi\Psi_j(x, \xi)\overline{\Phi_j(y, \xi)}
= \sum_{\xi\in \cX_j}\psi_\xi(x)\overline{\ph_\xi(y)}.
\end{eqnarray*}
Hence,
$$
\Psi_j*\overline{\Phi}_j*f = \sum_{\xi\in \cX_j}
\langle f,\ph_\xi\rangle \psi_\xi.
$$
Substituting this in (\ref{Needle-rep}) yields (\ref{needlet-rep1}).

(b) To prove (\ref{Needle-rep}) in $\Lpp$
we observe that
$\sum_{j=0}^{\ell} \Psi_j*\overline{\Phi}_j*f = \Ker_\ell*f$
with
$\Ker_\ell:=\sum_{j=0}^{\ell}\Psi_j\ast\overline{\Phi}_j$.
Because of (\ref{part-unity}), $\Ker_\ell(x, y)$
is a reproducing kernel for polynomials
exactly as the kernels $\Ker_n(x, y)$ from Lemma~\ref{lem:Ker-n}.
Hence, $\sum_{j=0}^\ell \Psi_j*\overline{\Phi}_j*f \to f$
in $\Lpp$ ($1\le p \le \infty$).
Then (\ref{needlet-rep1}) in $\Lpp$ follows as above.
The unconditional convergence in $\Lpp$, $1<p<\infty$, follows by
Proposition~\ref{t:identification} and
Theorem~\ref{thm:Fnorm-equivalence} below.
\qed


\begin{remark}\label{frame}
It is easy to see that there exists a function $\ha \ge 0$ satisfying
$(\ref{ha-hb1})-(\ref{ha-hb2})$ such that
$\ha^2(t) + \ha^2(2t) =1$, $t \in [1/2, 1]$.
Suppose that in the above construction
$\hb = \ha$ and $\ha\ge 0$.
Then
$\Phi_j = \Psi_j$ and
$\ph_\xi = \psi_\xi$.
Now $(\ref{needlet-rep1})$ becomes
$f = \sum_{\xi \in \cX} \langle f, \psi_\xi\rangle \psi_\xi$.
It is easily seen that
$\{\psi_\xi: \xi \in \cX\}$ is a tight frame for $\LL2$
$($see \cite{PX1}$)$.
\end{remark}

\section{First scale of weighted Triebel-Lizorkin  spaces on $[-1,1]$}\label{Tri-Liz}
\setcounter{equation}{0}

In analogy to the classical case on $\Rd$
we will define our first scale of weighted Triebel-Lizorkin spaces by means of
the Littlewood-Paley expressions employing the kernels $\Phi_j$,
defined by
\begin{equation}\label{def-Phi-j}
\Phi_0(x, y) := \PP_0(x)\PP_0(y)
\quad\mbox{and}\quad
\Phi_j(x, y) := \sum_{\nu=0}^\infty
\ha \Big(\frac{\nu}{2^{j-1}}\Big)\PP_\nu(x)\PP_\nu(y),
~~~j\ge 1,
\end{equation}
where $\ha$ satisfies the conditions
\begin{align}
&\quad \ha\in C^\infty[0, \infty),
\quad \supp \ha \subset [1/2, 2], \label{ha1}\\
&\quad  |\ha(t)|>c>0, \quad \text{if } t \in [3/5, 5/3].\label{ha2}
\end{align}

\begin{definition}
Let  $s \in \R$, $0<p<\infty$, and $0<q\le\infty$.
Then the weighted Triebel-Lizorkin space
$\F spq:=\F spq(\w)$ is defined as the set of all $f\in \cD'$ such that
\begin{equation}\label{Tri-Liz-norm}
\norm{f}_{\F spq}:=\nnorm{\biggl(\sum_{j=0}^{\infty}
(2^{sj}|\Phi_j\ast f(\cdot)|)^q\biggr)^{1/q}}_{\Lp} <\infty
\end{equation}
with the usual modification when $q=\infty$.
\end{definition}

Observe that the above definition is independent of the choice
of $\ha$ as long as it satisfies
(\ref{ha1})-(\ref{ha2}) (see Theorem~\ref{thm:Fnorm-equivalence} below).


\begin{proposition}\label{Fspq_embedding}
For every $s \in \R$, $0<p<\infty$, and $0<q\le\infty$,
$\F spq$ is a quasi-Banach space which is continuously embedded in $\cD'$.
\end{proposition}

\noindent
{\bf Proof.}
We will only prove the continuous embedding of $\F spq$ in $\cD'$.
Then the completeness follows by a standard argument (see e.g. \cite{T1}, p. 49).

Suppose the kernels $\Phi_j$ are as in the definition of $\F spq$
with $\ha$ satisfying (\ref{ha1})-(\ref{ha2}) which are the same as
(\ref{ha-hb1})-(\ref{ha-hb2}).
Then as was already mentioned, there is a function $\hb$ satisfying
(\ref{ha-hb1})-(\ref{ha-hb2}) such that (\ref{ha-hb3}) holds as well.
Let $\Psi_j$ be defined by (\ref{def-Psi-j}).
Then by Proposition~\ref{prop:needlet-rep} $f\in\F spq$ has the representation:
$
f=\sum_{j=0}^\infty \Psi_j*\Phi_j*f
$
in $\cD'$.
Hence for $\phi\in\cD$ we have
$
\langle f, \phi\rangle=\sum_{j=0}^\infty \langle \Psi_j*\overline{\Phi}_j*f, \phi\rangle.
$
Using (\ref{conv})-(\ref{conv-conv}) it follows that
$$
\langle \Psi_j\ast\overline{\Phi}_j\ast f, \phi\rangle
=\sum_{\nu=0}^{2^j} \overline{\ha\Big(\frac{\nu}{2^{j-1}}\Big)}
\hb\Big(\frac{\nu}{2^{j-1}}\Big)a_\nu(f)a_\nu(\phi)
$$
and
\begin{equation}
\begin{aligned}\label{embed9}
|\langle \Psi_j\ast\overline{\Phi}_j\ast f, \phi\rangle|
&\le c 2^{j/2}\|\Phi_j\ast f\|_\L2\max_{0\le \nu\le 2^j}|a_\nu(\phi)|\\
&\le c2^{j(2/p+1/2)}\|\Phi_j\ast f\|_\Lp\max_{0\le \nu\le 2^j}|a_\nu(\phi)|\\
&\le c2^{-j}\|f\|_{\F spq} \cN_k(\phi),
\end{aligned}
\end{equation}
where $k\ge 2/p+3/2-s$, $\cN_k(\cdot)$ is from (\ref{D-norms}),
and we used inequality (\ref{norm-relation}).
Consequently,
$|\langle f, \phi\rangle| \le c\|f\|_{\F spq}\cN_k(\phi)$,
which is the desired embedding.
$\qed$

\smallskip

It is natural to define the weighted potential space
(generalized weighted Sobolev space)
$H^p_s := H^p_s(\w)$,
$s >0$, $1\le p \le \infty$,  on $[-1, 1]$
as the set of all $f\in \cD'$ such that
\begin{equation}\label{def-H}
\|f\|_{H^p_s}:=
\Big\|\sum_{n=0}^\infty (n+1)^s a_n(f)\PP_n(\cdot) \Big\|_\Lp
< \infty,
\end{equation}
where
$a_n(f):= \langle f, \PP_n\rangle$ as in Lemma~\ref{lem:dec-D1}.

In the next statement we identify certain weighted Triebel-Lizorkin spaces
as weighted potential spaces or $\Lpp$.


\begin{proposition}\label{t:identification}
We have
\begin{equation}\label{ident1}
F^{s 2}_p \sim H^p_s,
\quad s >0, ~1 < p < \infty,
\end{equation}
and
\begin{equation}\label{ident2}
F^{0 2}_p \sim H^p_0 \sim \Lpp,
\quad 1 < p < \infty,
\end{equation}
with equivalent norms.
\end{proposition}

One proves this proposition in a standard way using e.g.
the multipliers from~\cite{CS}.
The proof can be carried out exactly as in the case of spherical
harmonic expansions, given in \cite[Proposition~4.3]{NPW},
and will be omitted.

Associated with $\F spq$ is the sequence space $\f spq$ defined as follows.

\begin{definition}
Let $s \in \R$, $0<p<\infty$, and $0<q\le\infty$. Then $\f spq$
is defined as the space of all complex-valued sequences
$h:=\{h_{\xi}\}_{\xi\in \cX}$ such that
\begin{equation}\label{def-f-space}
\norm{h}_{\f spq}
:=\nnorm{\biggl(\sum_{j=0}^\infty 2^{sjq}\sum_{\xi \in \cX_j}
(|h_{\xi}|\tONE_{I_\xi}(\cdot))^q\biggr)^{1/q}}_{\Lp} <\infty
\end{equation}
with the usual modification for $q=\infty$.
Recall that
$\tONE_{I_\xi}:=\mu(I_\xi)^{-1/2}\ONE_{I_\xi}$.
\end{definition}

We now introduce the ``analysis" and ``synthesis" operators
\begin{equation}\label{anal_synth_oprts}
S_\varphi: f\rightarrow \{\ip{f, \varphi_\xi}\}_{\xi \in \cX}
\quad\text{and}\quad
T_\psi: \{h_\xi\}_{\xi \in \cX}\rightarrow \sum_{\xi\in \X}h_\xi\psi_\xi.
\end{equation}


Here is our main result concerning the weighted $F$-spaces.


\begin{theorem}\label{thm:Fnorm-equivalence}
Let $s\in \R$, $0< p< \infty$ and $0<q\le \infty$. Then the operators
$S_\varphi:\F spq\rightarrow\f spq$ and $T_\psi:\f spq\rightarrow \F spq$
are bounded and $T_\psi\circ S_\varphi=Id$ on $\F spq$.
Consequently, for $f\in \cD'$ we have that $f\in \F spq$ if and only if
$\{\ip{f, \varphi_\xi}\}_{\xi \in \cX}\in \f spq$.
Furthermore,
\begin{align}\label{Fnorm-equivalence-1}
\norm{f}_{\F spq} &
\sim  \norm{\{\ip{f,\varphi_\xi}\}}_{\f spq}
\sim \nnorm{\biggl(\sum_{j=0}^\infty 2^{sjq}
\sum_{\xi\in \cX_j}|\ip{f, \varphi_\xi}\psi_\xi(\cdot)|^q\biggr)^{1/q}}_{\Lp}.
\end{align}
In addition, the definition of $\F spq$ is independent of the selection of
$\ha$ satisfying $(\ref{ha1})$--$(\ref{ha2})$.
\end{theorem}

For the proof of this theorem we will need several lemmas whose proofs are given
in the appendix.


\begin{lemma}\label{lem:Max-needl} If $\xi\in \cX_j,$ $j\ge 0$, and $0<t<\infty$,
then
\begin{equation}\label{Max-needl1}
|\ph_\xi(x)|, |\psi_\xi(x)|
\le c (\cM_t\tONE_{I_\xi})(x), \quad x\in [-1, 1], \quad\mbox{and}
\end{equation}
\begin{equation}\label{Max-needl2}
\tONE_{I_\xi}(x)
\le c(\cM_t \ph_\xi)(x), c(\cM_t \psi_\xi)(x), \quad x\in [-1, 1].
\end{equation}
\end{lemma}


\begin{lemma}\label{Almdiag}
For any $\sigma > 0$ there exists a constant $c_\sigma>0$ such that
\begin{equation}\label{almdiag}
|\Phi_j\ast \psi_\xi(x)|
\le c_\sigma \frac{2^{j/2}}{\sqrt{ \W(2^j;x) }(1+2^jd(\xi, x))^\sigma},
\quad \xi\in\cX_\nu,
\quad  j-1\le \nu\le j+1,
\end{equation}
and $\Phi_j\ast \psi_\xi(x)=0$ for $\xi\in\cX_\nu$,
$\nu\ge j+2$ or $\nu\le j-2$.
Here $\cX_\nu:=\emptyset$ if $\nu < 0$.
\end{lemma}


\begin{definition}
For a collection of complex numbers
$\{h_{\xi}\}_{\xi\in \cX_j}$ we let
\begin{equation}\label{def.h-star}
h^{\ast}_\xi:=\sum_{\eta\in
\cX_j}\frac{|h_\eta|}{(1+2^jd(\eta,\xi))^\sigma}.
\end{equation}
Here $\sigma>1$ is sufficiently large and will be selected later on.
\end{definition}


\begin{lemma}\label{l:weak_inequality}
Suppose that $P\in \Pi_{2^j},j\ge 0$, and
let $a_\xi:=\max_{x\in I_\xi}|P(x)|$ .
There exists $r\ge 1$, depending only on $\sigma$, $\a$, and $\b$,
such that if
$$
b_\xi:=\max\{\min_{x\in I_\eta }|P(x)|:\eta\in \cX_{j+r},
I_\xi\cap I_\eta\ne\emptyset \},
$$
then
\begin{equation}\label{weak_inequality}
a_\xi^\ast\sim b_\xi^\ast
\end{equation}
with constants of equivalence independent of $P,j$ and $\xi$.
\end{lemma}


\begin{lemma}\label{lem:disc-max}
Assume $t>0$ and let $\{b_\xi\}_{\xi\in \cX_j}$ $(j\ge 0)$
be a collection of complex numbers.
Suppose that
$\sigma>(4\max\{\a, \b\}+3)/t+1$
in the definition $(\ref{def.h-star})$ of $b_\xi^*$.
Then
\begin{equation*}
b_\xi^{\ast}\ONE_{I_\xi}(x)
\le c \cM_t\Big(\sum_{\eta\in \cX_j}|b_\eta|\ONE_{I_\eta}\Big)(x),
\quad x\in I_\xi, \quad \xi \in \cX_j.
\end{equation*}
\end{lemma}


\noindent
{\bf Proof of Theorem \ref{thm:Fnorm-equivalence}.}
Suppose $\a\ge\b$.
Fix  $0<t<\min\{p,q\}$ and let $\sigma>(4\a+3)/t+1$.
We first note that the right-hand side equivalence in (\ref{Fnorm-equivalence-1})
follows immediately from Lemma \ref{lem:Max-needl}
and the maximal inequality (\ref{max_ineq}).

Assume that $\{\Phi_j\}$ are from the definition of
weighted Triebel-Lizorkin spaces, i.e. $\Phi_j$ are defined by (\ref{def-Phi-j}),
where $\ha$ satisfies (\ref{ha1})-(\ref{ha2}),
the same as (\ref{ha-hb1})-(\ref{ha-hb2}).
As already mentioned, there exists a function $\hb$ satisfying
(\ref{ha-hb1})-(\ref{ha-hb2}) such that (\ref{ha-hb3}) holds.
Let $\Psi_j$ be defined by (\ref{def-Psi-j}) using this $\hb$.
Also, let $\{\ph_\xi\}_{\xi\in\cX}$ and $\{\psi_\xi\}_{\xi\in\cX}$
be the associated needlet systems defined as in (\ref{def-needlets1}).

Further, let $\{\wt\Phi_j\}$ be a second sequence of kernels
like the kernels $\{\Phi_j\}$ from above but defined by a different function $\ha$.
Also, we assume that a sequence of companion kernels $\{\wt\Psi_j\}$ is constructed
as above and let
$\{\wt\ph_\xi\}$, $\{\wt\psi_\xi\}$
be the associated needlet systems, defined as in
(\ref{def.Phi-j})-(\ref{def-needlets1}).
So, we have two totaly different systems of kernels and associated needlet systems.

We first establish the boundedness of $T_{\wt\psi}:\f spq\rightarrow \F spq$,
defined similarly as in~(\ref{anal_synth_oprts}),
where the space $\F spq$ is defined by $\{\Phi_j\}$.
Let $h:=\{h_\xi\}_{\xi\in \cX}$
be a finitely supported sequence and  $f:=\sum_{\xi}h_\xi \psi_\xi$.
Using (\ref{almdiag}) we have for $x\in [-1,1]$,
\begin{align*}
|\Phi_j\ast f(x)|
&=\Big|\sum_{\xi\in \cX}  h_\xi \Phi_j\ast\wt\psi_{\xi}(x)\Big|
\le \sum_{j-1\le \nu\le j+1}\sum_{\xi\in \cX_\mu} |h_\xi |
|\Phi_j\ast\wt\psi_{\xi}(x)|\\
&\le c2^{j/2} \sum_{j-1\le \nu\le j+1}\sum_{\xi\in \cX_\nu}
\frac{ |h_\xi| }{\sqrt{ \W(2^\nu;x)}(1+2^\nu d(\xi, x))^\sigma}.
\end{align*}
Fix $\eta\in\cX_j$ and denote
$\cY_\eta:=\{\xi\in\cX_{j-1}\cup\cX_j\cup\cX_{j+1}: I_\xi\cap I_\eta\ne \emptyset\}$
($\cX_{-1}:=\emptyset$).
Notice that $\# \cY_\eta \le {\rm constant}$ and
$d(x, \xi) \le c2^{-j}$ if $x\in I_\xi$ and $\xi\in\cY_\eta$.
Hence, we have for $x\in I_\eta$
\begin{align*}
|\Phi_j\ast f(x)|
&\le c2^{j/2} \sum_{j-1\le \nu\le j+1}
\sum_{\omega\in\cY_\eta\cap\cX_\nu}\sum_{\xi\in \cX_\nu}
\frac{|h_\xi|\ONE_{\omega}(x)}{\sqrt{\W(2^\nu;\omega)}(1+2^\nu d(\xi, \omega))^\sigma}\\
&\le c2^{j/2}
\sum_{\omega\in\cY_\eta}\frac{h_\omega^*\ONE_{\omega}(x)}{\sqrt{ \W(2^j;\omega)}}
\le c\sum_{\omega\in\cY_\eta}h_\omega^*\tONE_{\omega}(x),
\end{align*}
where we also used (\ref{muI}).
We now insert this in (\ref{Tri-Liz-norm})
and use Lemma \ref{lem:disc-max} and the maximal inequality (\ref{max_ineq})
to obtain
\begin{equation}
\begin{aligned}\label{F-less-f}
\norm{f}_{\F spq}
&\le c\nnorm{\Bigr(\sum_{j=0}^\infty \Big[2^{sj}\sum_{\eta\in \cX_j}\sum_{\omega\in \cY_\eta}
h_\omega^\ast\tONE_{I_\omega}(\cdot)\Big]^q\Big)^{1/q}}_{\Lp}\\
&\le c\nnorm{\Bigr(\sum_{j=0}^\infty \Big[2^{sj}\sum_{\xi\in \cX_j}
h_\xi^\ast\tONE_{I_\xi}(\cdot)\Big]^q\Big)^{1/q}}_{\Lp}\\
&\le c\nnorm{\Bigr(\sum_{j=0}^\infty \Big[\cM_t\Big(\sum_{\xi\in \cX_j}
2^{sj}|h_\xi|\tONE_{I_\xi}\Big)(\cdot)\Big]^q\Big)^{1/q}}_{\Lp}\\
&\le c\norm{\{h_\xi\} }_{\f spq}.
\end{aligned}
\end{equation}
For the second estimate above it was important that $\#\cY_\eta\le c$.
This establishes the desired result for finitely supported sequences.
Using  the continuous  embedding of $\F spq$ in $\cD'$
(Lemma \ref{Fspq_embedding}) and the density of finitely supported
sequences in $\f spq$ it follows from  (\ref{F-less-f})
that for every $h\in \f spq$,
$T_{\wt\psi} h:=\sum_{\xi\in \cX}h_\xi\wt\psi_\xi$
is a well defined distribution in $\cD'$.
Thus a standard density argument
shows that $T_{\wt\psi}:\f spq\rightarrow\F spq$ is bounded.


We next prove the boundedness of the operator $S_\varphi:\F spq\rightarrow \f spq$,
where we assume this time that $\F spq$ is defined in terms
of $\{\overline{\Phi}_j\}$.
Let $f\in \F spq$. Then $\overline{\Phi}_j\ast f\in \Pi_{2^j}$.
For $\xi \in \cX_j$, we set
$$
a_\xi:=\max_{x\in
I_\xi}|\overline{\Phi}_j\ast f(x)|, \quad
b_\xi:=\max\{\min_{x\in I_\eta }|\overline{\Phi}_j\ast f(x)|:\eta\in \cX_{j+r},
I_\xi\cap I_\eta\ne\emptyset \}.
$$
Assuming that $r$ above is the constant from Lemma \ref{l:weak_inequality},
it follows by the same lemma that $a_\xi^*\sim b_\xi^*$.
Therefore,
$$
|\langle f, \ph_\xi \rangle|=\cc_\xi^{1/2}|\overline{\Phi}_j*f(\xi)|
\le c\mu(I_\xi)^{1/2}a_\xi \le c\mu(I_\xi)^{1/2}a_\xi^*
\le c\mu(I_\xi)^{1/2}b_\xi^*.
$$
>From this, taking into account that $\tONE_{I_\xi}:=\mu(I_\xi)^{-1/2}\ONE_{I_\xi}$,
we obtain
\begin{equation}\label{asdfasf}
\begin{aligned}
\norm{\{\ip{f, \varphi_\xi}\}}_{\f spq}
&:=\nnorm{\Bigl(\sum_{j\ge 0} 2^{jsq}
\sum_{\xi\in \cX_j}
[|\ip{f, \varphi_\xi}|\tONE_{I_\xi}(\cdot)]^q\Bigr)^{1/q}}_{\Lp}\\
&\le c\nnorm{\Bigl(\sum_{j\ge 0} 2^{jsq}
\sum_{\xi\in \cX_j}[b^\ast_\xi
\ONE_{I_\xi}(\cdot)]^q\Bigr)^{1/q}}_{\Lp}\\
&\le c\nnorm{\Bigl(\sum_{j\ge 0}2^{jsq}
\Big[\cM_t\Big(\sum_{\xi\in \cX_j}
b_\xi\ONE_{I_\xi}\Big)(\cdot)\Big]^q\Bigr)^{1/q}}_{\Lp}\\
&\le c\nnorm{\Bigl(\sum_{j\ge 0}2^{jsq}
\sum_{\xi\in \cX_j} b_\xi^q\ONE_{I_\xi}(\cdot)\Bigr)^{1/q}}_{\Lp},
\end{aligned}
\end{equation}
where for the second inequality above we used Lemma~\ref{lem:disc-max}
and for the third the maximal inequality (\ref{max_ineq}).

Let
$m_\eta:=\min_{x\in I_\eta}|\overline{\Phi}_j\ast f(x)|$ for $\xi\in \cX_{j+r}$
and denote, for $\xi\in\cX_j$,
$$\cX_{j+r}(\xi):=\{w\in \cX_{j+r}:I_w\cap I_\xi \ne \emptyset\}.$$
Evidently, $\#\cX_{j+r}(\xi)\le \wt c$, $\wt c= \wt c(r)$.
Hence, $d(w,\eta)\le c(r)2^{-j-r}$
for $w, \eta\in \cX_{j+r}(\xi)$ and therefore
$$
m_w\le c\frac{m_w}{1+2^{j+r}d(w,\eta)}\le cm^\ast_\eta.
$$
Consequently, for any $\xi\in \cX_j$ and $\eta\in \cX_{j+r}(\xi)$,
we have $b_\xi=\max_{w\in\cX_{j+r}(\xi)}m_w \le cm^\ast_\eta$
and hence
$$
b_\xi\ONE_{I_\xi}\le c\sum_{\eta\in \cX_{j+r}(\xi)}  m^\ast_\eta \ONE_{I_\eta}.
$$
Using  this estimate in (\ref{asdfasf}) we get
\begin{equation*}
\begin{aligned}
\norm{\{\ip{f, \varphi_\xi}\}}_{\f spq}
%
&\le c\nnorm{\Bigl(\sum_{j\ge 0}2^{jsq}
\Bigl(\sum_{\eta\in \cX_{j+r}}  m^\ast_\eta
\ONE_{I_\eta}(\cdot) \Bigr)^q\Bigr)^{1/q}}_{\Lp}\\
&\le c\nnorm{\Bigl(\sum_{j\ge 0}2^{jsq}
\Big[\cM_t\Big(\sum_{\eta\in \cX_{j+r}}m_\eta
\ONE_{I_\eta}\Big)(\cdot)\Big]^q\Bigr)^{1/q}}_{\Lp}\\
&\le c\nnorm{\Bigl(\sum_{j\ge 0}\Bigl(2^{js}
\sum_{\xi\in \cX_j} m_\xi\ONE_{I_\xi}(\cdot) \Bigr)^q\Bigr)^{1/q}}_{\Lp}\\
&\le c\nnorm{\Bigl(\sum_{j\ge 0} (2^{js}
  |\overline{\Phi}_j\ast f|)^q\Bigr)^{1/q}}_{\Lp}.
\end{aligned}
\end{equation*}
Thus the boundedness of $S_\varphi:\F spq\rightarrow \f spq$ is established.

The identity $T_\psi\circ S_\varphi=Id$ follows by Proposition~\ref{prop:needlet-rep}.

It remains to show that $\F spq$ is  independent of the particular selection
of $\ha$ in the definition of  $\{\Phi_j\}$.
Denote for an instant by $\|f\|_{\F spq(\Phi)}$ the F-norm defined by $\{\Phi_j\}$.
Then by the above proof it follows that
$$
\|f\|_{\F spq(\Phi)}
\le c\|\{\langle f, \wt\ph_\xi\rangle\}\|_{\f spq}
\quad\mbox{and}\quad
\|\{\langle f, \ph_\xi\rangle\}\|_{\f spq}
\le c\|f\|_{\F spq(\overline{\Phi})}
$$
and hence
$$
\|f\|_{\F spq(\Phi)}
\le c\|\{\langle f, \wt\ph_\xi\rangle\}\|_{\f spq}
\le c\|f\|_{\F spq(\overline{\wt\Phi})}.
$$
Now the desired independence follows by
reversing the roles of $\{ \Phi_j\}$,$\{\wt\Phi_j\}$,
and  their complex conjugates.
$\qed$

\section{Second scale of weighted Triebel-Lizorkin  spaces on $[-1,1]$}\label{Tri-Liz-2}
\setcounter{equation}{0}

We introduce our second scale of Triebel-Lizorkin  spaces by
utilizing again the kernels
$\Phi_j$ defined by (\ref{def-Phi-j}) with $\ha$ satisfying
(\ref{ha1})-(\ref{ha2}) (compare with \S\ref{Tri-Liz}).

\begin{definition}
Let  $s \in \R$, $0<p<\infty$, and $0<q\le\infty$.
Then the weighted Triebel-Lizorkin space
$\FF spq:=\FF spq(\w)$ is defined as the set of all $f\in \cD'$ such that
\begin{equation}\label{Tri-Liz-norm2}
\norm{f}_{\FF spq}:=\nnorm{\biggl(\sum_{j=0}^{\infty}
\Big[2^{sj}\W(2^j;\cdot)^{-s}|\Phi_j\ast f(\cdot)|\Big]^q\biggr)^{1/q}}_{\Lp} <\infty
\end{equation}
with the usual modification when $q=\infty$.
\end{definition}

Observe that the above definition is independent of the choice
of $\ha$ as long as it satisfies
(\ref{ha1})-(\ref{ha2}) (see Theorem~\ref{thm:Fnorm-eq2} below).
Following in the footsteps of the development from \S\ref{Tri-Liz},
it is easy to show that $\FF spq$ is a complete quasi-Banach space,
which is embedded continuously in $\cD'$.
For the latter one proceeds as in the proof of Proposition~\ref{Fspq_embedding},
where in (\ref{embed9}) one, in addition, uses the obvious estimate
$\|g\|_2 \le cn^{\gamma}\|\W(n;\cdot)^{s} g(\cdot)\|_2$,
where $\gamma := (2\min\{\a,\b\}+1)s_+$,
which is immediate from
$c_1n^{-2\min\{\a,\b\}-1}\le \W(n; x)\le c_2$, $x\in [-1, 1]$.
We skip the details.

The sequence space $\ff spq$ associated with $\FF spq$ is now defined as follows.

\begin{definition}
Let $s \in \R$, $0<p<\infty$, and $0<q\le\infty$. Then $\ff spq$
is defined as the space of all complex-valued sequences
$h:=\{h_{\xi}\}_{\xi\in \cX}$ such that
\begin{equation}\label{def-f-space2}
\norm{h}_{\ff spq}
:=\nnorm{\biggl(\sum_{\xi \in \cX}
\Big[\mu(I_\xi)^{-s}|h_{\xi}|\tONE_{I_\xi}(\cdot)\Big]^q\biggr)^{1/q}}_{\Lp} <\infty
\end{equation}
with the usual modification when $q=\infty$.
\end{definition}

To characterize the Triebel-Lizorkin spaces $\FF spq$ we use again
the operators $S_\varphi$ and $T_\psi$ defined in (\ref{anal_synth_oprts})
and the sequence spaces $\ff spq$.


\begin{theorem}\label{thm:Fnorm-eq2}
Let $s\in \R$, $0< p< \infty$ and $0<q\le \infty$. Then the operators
$S_\varphi:\FF spq\rightarrow\ff spq$ and $T_\psi:\ff spq\rightarrow \FF spq$
are bounded and $T_\psi\circ S_\varphi=Id$ on $\FF spq$.
Consequently, for $f\in \cD'$ we have that $f\in \FF spq$ if and only if
$\{\ip{f, \varphi_\xi}\}_{\xi \in \cX}\in \ff spq$.
Furthermore,
\begin{align}\label{Fnorm-equivalence-12}
\norm{f}_{\FF spq} &
\sim  \norm{\{\ip{f,\varphi_\xi}\}}_{\ff spq}
\sim \nnorm{\biggl(\sum_{\xi\in \cX}
\Big[\mu(I_\xi)^{-s}|\ip{f, \varphi_\xi}\psi_\xi(\cdot)|\Big]^q\biggr)^{1/q}}_{\Lp}.
\end{align}
In addition, the definition of $\FF spq$ is independent of the selection of
$\ha$ satisfying $(\ref{ha1})$--$(\ref{ha2})$.
\end{theorem}

The proof of this theorem is similar to the proof of
Theorem~\ref{thm:Fnorm-equivalence}.
The only new ingredient is the following lemma.


\begin{lemma}\label{lem:disc-max2}
Let $t>0$ and $s\in \RR$. Suppose $\{b_\xi\}_{\xi\in \cX_j}$ $(j\ge 0)$
is a collection of complex numbers and let
$\sigma>(4\max\{\a, \b\}+3)(1/t+|s|)+1$
in the definition $(\ref{def.h-star})$ of $b_\xi^*$.
Then
\begin{equation}\label{disc-max22}
\mu(I_\xi)^{-s}b_\xi^{\ast}\ONE_{I_\xi}(x)
\le c \cM_t\Big(\sum_{\eta\in \cX_j}\mu(I_\eta)^{-s}|b_\eta|\ONE_{I_\eta}\Big)(x),
\quad x\in I_\xi, \quad \xi \in \cX_j,
\end{equation}
\end{lemma}

\noindent
{\bf Proof.}
For $\xi\in\cX_j$, $\mu(I_\xi)\sim 2^{-j}\W(2^j;\xi)$ and hence,
using (\ref{omega<omega}),
\begin{align*}
\mu(I_\xi)^{-s}b_\xi^*
\le c\sum_{\eta\in\cX_j}\frac{2^{js}\W(2^j;\xi)^{-s}|b_\eta|}
                              {(1+2^jd(\xi, \eta))^\sigma}
\le c\sum_{\eta\in\cX_j}\frac{2^{js}\W(2^j;\eta)^{-s}|b_\eta|}
                              {(1+2^jd(\xi, \eta))^{\sigma_1}}
\le c\Big(\mu(I_\eta)^{-s}|b_\eta|\Big)^*,
\end{align*}
where $\sigma_1:=\sigma-(2\max\{\a,\b\}+1)|s|>(4\max\{\a, \b\}+3)/t+1$.
Now (\ref{disc-max22}) follows by Lemma~\ref{lem:disc-max}.
$\qed$

Now the proof of Theorem~\ref{thm:Fnorm-eq2} can be carried out as the proof of
Theorem~\ref{thm:Fnorm-equivalence}, using Lemma~\ref{lem:disc-max2}
in place of Lemma~\ref{lem:disc-max} and selecting $\sigma$ in
the definitions of $h_\xi^*$ and $a_\xi^*$, $b_\xi^*$ sufficiently large.
We skip the further details.

\smallskip

In a sense the spaces $\FF spq$ are more natural than the spaces $\F spq$
from \S\ref{Tri-Liz}
since they scale (are embedded) ``correctly" with respect to the smoothness
index $s$.

%
\begin{proposition}\label{F-embedding}
Let $0<p<p_1<\infty$, $0<q, q_1\le\infty$, and $0<s_1<s<\infty$.
Then we have the continuous embedding
\begin{equation}\label{F-embed}
\FF spq \subset \FF {s_1}{p_1}{q_1}
\quad\mbox{if}\quad s-1/p=s_1-1/p_1.
\end{equation}
\end{proposition}
The proof of this embedding result can be carried out similarly as
in the classical case on $\RR^n$ using inequality (\ref{norm-relation2})
and Theorem~\ref{thm:Fnorm-eq2} (see e.g. \cite{T1}, page 129).
It~will be omitted.

\section{First scale of weighted Besov spaces on $[-1,1]$}\label{Besov}
\setcounter{equation}{0}

To introduce the first scale of weighted Besov spaces we use the kernels
$\Phi_j$ defined in (\ref{def-Phi-j}) with $\ha$ satisfying
(\ref{ha1})-(\ref{ha2}) (see \cite{RS, T1}).

\begin{definition}
Let $s\in \R$ and $0<p,q \le \infty$. Then the weighted Besov space
$\Bes spq := \Bes spq(\w)$ is defined
as the set of all $f \in \cD'$ such that
\begin{equation*}
\|f\|_{\Bes spq} :=
\Big(\sum_{j=0}^\infty \Big(2^{s j}
\|\Phi_j\ast f\|_{\Lp}\Big)^q\Big)^{1/q}
< \infty,
\end{equation*}
where the $\ell_q$-norm is replaced by the sup-norm if $q=\infty$.
\end{definition}
Note that as in the case of weighted Triebel-Lizorkin spaces the above definition
is independent of the choice of $\ha$ satisfying  (\ref{ha1})-(\ref{ha2})
(see Theorem~\ref{thm:Bnorm-eq}).
Also, the Besov space $\Bes spq(\w)$ is a quasi-Banach
space which is continuously embedded in $\cD'$.


Our next goal is to link the weighted Besov spaces
with best polynomial approximation in $\Lpp$.
Recall that $E_n(f)_p$ denotes the best approximation of $f \in \Lpp$
from $\Pi_n$ (see (\ref{def-En})).


\begin{proposition}\label{prop:character-Besov}
Let $s > 0$, $1 \le p \le \infty$, and $0 < q \le \infty$.
Then $f \in B^{s q}_p$ if and only if
\begin{equation}\label{character-Besov1}
\|f\|_{B^{s q}_p}^A
:= \|f\|_\Lp +
\Big(\sum_{j=0}^\infty (2^{s j}E_{2^j}(f)_p)^q\Big)^{1/q}
< \infty.
\end{equation}
Moreover,
\begin{equation}\label{character-Besov2}
\|f\|_{B^{s q}_p}^A \sim \|f\|_{B^{s q}_p}.
\end{equation}
\end{proposition}

\noindent
{\bf Proof.}
Let $f \in B^{s q}_p$.
It is easy and standard to show that under the assumptions on
$s$, $p$, and $q$ the space
$B^{s q}_p$ is continuously imbedded in $\Lpp$,
i.e. $f$ can be identified as a function in $\Lpp$ and
$\|f\|_\Lp \le c\|f\|_{B^{s q}_p}$.

It is easy to construct (see e.g. \cite{FJ1}) a function $\ha\ge 0$
satisfying (\ref{ha1})-(\ref{ha2}) such that
$\ha(t) + \ha(2t)  =1$ for $t \in [1/2, 1]$ and hence
\begin{equation}\label{unity0}
\sum_{\nu=0}^\infty \ha(2^{-\nu}t) =1,
\quad t \in [1, \infty).
\end{equation}
Assume that $\{\Phi_j\}$ are defined by
(\ref{def-Phi-j}) with such a function $\ha$.
As in Proposition~\ref{prop:needlet-rep}, it is easy to see that
$f = \sum_{j=0}^\infty \Phi_j*f$ in $\Lpp$.
Hence, since $\Phi_j*f \in \Pi_{2^j}$,
\begin{equation}\label{E2j<sum}
E_{2^\ell}(f)_p \le \sum_{j=\ell+1}^\infty \|\Phi_j*f\|_p,
\quad \ell\ge 0.
\end{equation}
Now, a standard argument using (\ref{E2j<sum}) shows that
$\|f\|_{B^{s q}_p}^A \le c\|f\|_{B^{s q}_p}$.

To prove the estimate in the other direction, we note that
$\Phi_j*f = \Phi_j*(f-Q)$ for $Q \in \Pi_{2^{j-2}}$ ($j\ge 2$).
Hence, as in Lemma~\ref{lem:Ker-n},
$
\|\Phi_j*f\|_\Lp \le c\|f-Q\|_\Lp.
$
Therefore,
$$
\|\Phi_j*f\|_\Lp \le cE_{2^{j-2}}(f)_p,
\quad j\ge 2,
\quad\mbox{and}\quad
\|\Phi_j*f\|_\Lp \le c\|f\|_\Lp,
$$
which imply
$
\|f\|_{B^{s q}_p} \le c\|f\|_{B^{s q}_p}^A.
$

Above we used that the definition of $B^{s q}_p$ is independent of
the selection of $\ha$, satisfying (\ref{ha1})-(\ref{ha2}).
$\qed$


\begin{remark}\label{moduli}
Proposition~\ref{prop:character-Besov} shows that when
$s > 0$ and $1 \le p \le \infty$
the weighted Besov spaces $B^{s q}_p$ can be identified as approximation
spaces induced by best polynomial approximation in $\Lpp$.
Also, it is worth mentioning that $E_n(f)_p$ can be characterized via
the weighted moduli of smoothness of Ditzian-Totik \cite{DT}.
Consequently, the weighted moduli of smoothness can be used for characterization
of weighted Besov spaces as well.
\end{remark}


It is natural to associate with the weighted Besov space $\Bes spq$ the
sequence space $\bes spq$ defined as follows.
\begin{definition}
Let $s\in \R$ and $0<p,q \le \infty$. Then $\bes spq:=\bes spq(\w)$
is defined to be the space of all complex-valued sequences
$h:=\{h_{\xi}\}_{\xi\in \cX}$ such that
\begin{equation*}
\norm{h}_{\bes spq}
:=\Bigl(\sum_{j=0}^\infty 2^{jsq}\Bigl[\sum_{\xi\in \cX_j}
\Big(\mu(I_\xi)^{1/p-1/2}|h_\xi|\Big)^p \Bigr]^{q/p}\Bigr)^{1/q}<\infty
\end{equation*}
with the usual modification for $p=\infty$ or $q=\infty$.
\end{definition}

Our main result in this section is the following characterization of
weighted Besov spaces,
which employs the operators $S_\ph$ and $T_\psi$ defined in (\ref{anal_synth_oprts}).


\begin{theorem}\label{thm:Bnorm-eq}
Let $s\in \R$ and  $0< p,q\le \infty$.
The operators $S_\varphi:\Bes spq\rightarrow\bes spq$ and
$T_\psi:\bes spq\rightarrow \Bes spq$ are bounded and
$T_\psi\circ S_\varphi=Id$ on $\Bes spq$.
Consequently, for $f\in \cD'$ we have that $f\in \Bes spq$ if and only if
$\{\ip{f, \varphi_\xi}\}_{\xi \in \cX}\in \bes spq$.
Moreover,
\begin{align}\label{Bnorm-equivalence-1}
\norm{f}_{\Bes spq} &
\sim  \norm{\{\ip{f,\varphi_\xi}\}}_{\bes spq}
\sim \Big(\sum_{j=0}^\infty2^{sjq}
\Bigl[\sum_{\xi\in \cX_j}
\norm{\ip{f,\varphi_\xi}\psi_\xi}_{\Lp}^p\Bigr]^{q/p}\Bigr)^{1/q}.
\end{align}
In addition, the definition of $\Bes spq$ is independent of the selection of
$\ha$ satisfying $(\ref{ha1})$--$(\ref{ha2})$.
\end{theorem}

To prove this theorem we will need the following lemma whose proof
is presented in the appendix.


\begin{lemma}\label{l:half_shannon}
For every $P\in \Pi_{2^j}, j\ge 0$, and $0<p\le \infty$
\begin{equation}
 \Big(\sum_{\xi\in \cX_j}\max_{x\in I_\xi}|P(x)|^p \mu(I_\xi)\Big)^{1/p}
 \lesssim   \norm{P}_{\Lp}.
\end{equation}
\end{lemma}


\noindent
{\bf Proof of Theorem~\ref{thm:Bnorm-eq}.}
Note first that the right-hand side of (\ref{Bnorm-equivalence-1})
follows immediately from (\ref{norm-Needlets}).

As in the proof of Theorem~\ref{thm:Fnorm-equivalence},
assume that the kernels $\Phi_j$ are defined by (\ref{def-Phi-j}),
where $\ha$ satisfies (\ref{ha1})-(\ref{ha2}).
Let $\hb$ be such that (\ref{ha-hb1})-(\ref{ha-hb3}) hold
and let $\Psi_j$ be defined by (\ref{def-Psi-j}) using this $\hb$.
Also, let $\{\ph_\xi\}_{\xi\in\cX}$ and $\{\psi_\xi\}_{\xi\in\cX}$
be the associated needlet systems defined as in (\ref{def-needlets1}).
Further, assume that
$\{\wt\Phi_j\}$, $\{\wt\Psi_j\}$, $\{\wt\ph_\xi\}$, $\{\wt\psi_\xi\}$
is a second set of kernels and needlets.

We first prove the boundedness of the operator
$T_{\wt\psi}:\bes spq\rightarrow\Bes spq$,
where $\Bes spq$ is defined via $\{\Phi_j\}$.
Let $0<t<\min\{p,1\}$ and $\sigma \ge (2\a+2)/t+ \a +1/2$.
Assume that $h=\{h_\xi\}$ is a finitely supported sequence and
set $f:=\sum_{\xi\in \cX}h_\xi\psi_\xi$.
Employing Lemmata~\ref{lem:J-maximal}, \ref{Almdiag}, and
(\ref{omega<omega}) we get
\begin{align*}
|\Phi_j\ast f(x)|
%
&\le \sum_{j-1\le \nu\le j+1}\sum_{\xi\in \cX_\nu}
|h_\xi|| \Phi_\mu\ast\wt\psi_{\xi}(x)|\\
&\le c\sum_{j-1\le \nu\le j+1}\sum_{\xi\in \cX_\nu}
|h_\xi| \frac{2^{j/2} }{\sqrt{ \W(2^j;x)}(1+2^jd(\xi, x))^\sigma } \\
&\le c\sum_{j-1\le \nu\le j+1}\sum_{\xi\in \cX_\nu}
|h_\xi| \frac{2^{j/2} }{\sqrt{ \W(2^j;\xi)}(1+2^jd(\xi, x))^{\sigma-\a-1/2} } \\
%
&\le c\sum_{j-1\le \nu\le j+1}\sum_{\xi\in \cX_\nu}
|h_\xi| \mu(I_\xi)^{-1/2}\cM_t(\ONE_{I_\xi})(x),
\end{align*}
where we also used that $\sigma \ge (2\a+2)/t+ \a +1/2$.
Using the maximal inequality (\ref{max_ineq}) it follows that
\begin{equation*}
\begin{aligned}
\norm{\Phi_j\ast f}_{\Lp}^p
&\le \Big\|\sum_{j-1\le \nu\le j+1}\sum_{\xi\in \cX_\nu}
|h_\xi|\mu(I_\xi)^{-1/2}\cM_t(\ONE_{I_\xi})(\cdot)\Big\|_{\Lp}^p\\
&\le c\sum_{j-1\le \nu\le j+1}\sum_{\xi\in \cX_\nu}
|h_\xi|^p \mu(I_\xi)^{-p/2}\int_{-1}^1\ONE_{I_\xi}(x)\w(x)\, dx\\
&\le c\sum_{j-1\le \nu\le j+1}\sum_{\xi\in \cX_\nu}
|h_\xi| ^p  \mu(I_\xi)^{1-p/2}.
\end{aligned}
\end{equation*}
Multiplying by $2^{js}$ and summing over $j\ge0$ we get
$\norm{f}_{\Bes spq}\le c\norm{\{h_\xi\}}_{\bes spq}$.
To extend the result to an arbitrary sequence $h=\{h_\xi\}\in \bes spq$
one proceeds similarly as in the Triebel-Lizorkin case
by using the embedding of $\Bes spq$ in $\cD'$
and the density of finitely supported sequences in $\bes spq$.

We next prove the boundedness of the operator
$S_\varphi:\Bes spq\rightarrow \bes spq$,
where we assume that $\Bes spq$ is defined in terms
of $\{\overline{\Phi_j}\}$.
Note first that
$$
|\ip{f,\varphi_\xi}|\sim \mu(I_\xi)^{1/2}|\overline{\Phi_j} \ast f(\xi)|,
\quad \xi \in \cX_j.
$$
Since $\overline{\Phi_j}\ast f\in \Pi_{2^j}$, then using Lemma~\ref{l:half_shannon}
\begin{align*}
\sum_{\xi\in \cX_j}\mu(I_\xi)^{1-p/2}|\ip{f,\varphi_\xi}|^p
\le c\sum_{\xi\in \cX_j}\mu(I_\xi)\sup_{x\in I_\xi}|\overline{\Phi_j}\ast f(x)|^p
\le c\norm{\overline{\Phi_j}\ast f}_{\Lp}^p,
\end{align*}
which yields
$
\|\{\langle f, \ph\rangle\}\|_{\bes spq} \le c\norm{f}_{\Bes spq}.
$

The identity $T_\psi\circ S_\varphi=Id$ follows by Proposition~\ref{prop:needlet-rep}.

The independence of $\Bes spq$ of the particular selection
of $\ha$ in the definition of  $\{\Phi_j\}$ follows from above exactly as in
the Triebel-Lizorkin case (see the proof of Theorem~\ref{thm:Fnorm-equivalence}).
$\qed$

\section{Second scale of weighted Besov spaces on $[-1,1]$}\label{Besov2}
\setcounter{equation}{0}

We introduce a second scale of weighted Besov spaces
by using again as in \S\ref{Besov} the kernels
$\Phi_j$, defined by (\ref{def-Phi-j}) with $\ha$ satisfying
(\ref{ha1})-(\ref{ha2}).

\begin{definition}
Let $s\in \R$ and $0<p,q \le \infty$. Then the weighted Besov space
$\BBes spq := \BBes spq(\w)$ is defined
as the set of all $f \in \cD'$ such that
\begin{equation*}
\|f\|_{\BBes spq} :=
\Big(\sum_{j=0}^\infty \Big[2^{s j}
\|\W(2^j;\cdot)^{-s}\Phi_j*f(\cdot)\|_{\Lp}\Big]^q\Big)^{1/q}
< \infty,
\end{equation*}
where the $\ell_q$-norm is replaced by the sup-norm if $q=\infty$.
\end{definition}
As for the other weighted Besov and Triebel-Lizorkin spaces considered here
the above definition is independent of the choice of $\ha$ satisfying
(\ref{ha1})-(\ref{ha2}).
Also, the Besov space $\BBes spq(\w)$ is a quasi-Banach
space which is continuously embedded in $\cD'$.

The main advantages of the spaces $\BBes spq$ over $\Bes spq$ are that,
first, they scale (are embedded) ``correctly" with respect to the smoothness
index $s$, and secondly, the right smoothness spaces in nonlinear n-term
weighted approximation from needles are defined in terms of
spaces $\BBes spq$ (see \S\ref{Nonlin-app} below).

%
\begin{proposition}\label{B-embedding}
Let $0<p\le p_1<\infty$, $0<q\le q_1\le \infty$, and $0<s_1\le s<\infty$.
Then we have the continuous embedding
\begin{equation}\label{B-embed}
\BBes spq \subset \BBes {s_1}{p_1}{q_1}
\quad\mbox{if}\quad s-1/p=s_1-1/p_1.
\end{equation}
\end{proposition}
This embedding result follows readily by applying
inequality (\ref{norm-relation2}).


We now define the companion to $\BBes spq(\w)$
sequence space $\bbes spq(\w)$.
\begin{definition}
Let $s\in \R$ and $0<p,q \le \infty$. Then $\bbes spq:=\bbes spq(\w)$
is defined to be the space of all complex-valued sequences
$h:=\{h_{\xi}\}_{\xi\in \cX}$ such that
\begin{equation*}
\norm{h}_{\bbes spq}
:=\Bigl(\sum_{j=0}^\infty \Bigl[\sum_{\xi\in \cX_j}
\Big(\mu(I_\xi)^{-s+1/p-1/2}|h_\xi|\Big)^p \Bigr]^{q/p}\Bigr)^{1/q}<\infty
\end{equation*}
with the usual modification for $p=\infty$ or $q=\infty$.
\end{definition}

For the characterization of weighted Besov spaces $\BBes spq$,
we again employ the operators $S_\ph$ and $T_\psi$ defined in (\ref{anal_synth_oprts}).


\begin{theorem}\label{thm:Bnorm-eq2}
Let $s\in \R$ and  $0< p,q\le \infty$.
The operators $S_\varphi:\BBes spq\rightarrow\bbes spq$ and
$T_\psi:\bbes spq\rightarrow \BBes spq$ are bounded and
$T_\psi\circ S_\varphi=Id$ on $\BBes spq$.
Consequently, for $f\in \cD'$ we have that $f\in \BBes spq$ if and only if
$\{\ip{f, \varphi_\xi}\}_{\xi \in \cX}\in \bbes spq$.
Moreover,
\begin{align}\label{Bnorm-equivalence-12}
\norm{f}_{\BBes spq} &
\sim  \norm{\{\ip{f,\varphi_\xi}\}}_{\bbes spq}
\sim \Big(\sum_{j=0}^\infty
\Bigl[\sum_{\xi\in \cX_j}
\mu(I_\xi)^{-sp}\norm{\ip{f,\varphi_\xi}\psi_\xi}_{\Lp}^p\Bigr]^{q/p}\Bigr)^{1/q}.
\end{align}
In addition, the definition of $\BBes spq$ is independent of the selection of
$\ha$ satisfying $(\ref{ha1})$-$(\ref{ha2})$.
\end{theorem}

The following additional lemma is needed for the proof of
Theorem~\ref{thm:Bnorm-eq2}.


\begin{lemma}\label{lem:half-shannon2}
For every $P\in \Pi_{2^j}, j\ge 0$, and $0<p\le \infty$
\begin{equation}\label{half-shannon2}
\Big(\sum_{\xi\in \cX_j}\W(2^j;\xi)^{-sp}\sup_{x\in I_\xi}|P(x)|^p \mu(I_\xi)\Big)^{1/p}
\le c\|\W(2^j;\cdot)P(\cdot)\|_{\Lp}.
\end{equation}
\end{lemma}

The proof of this lemma is similar to the proof of Lemma~\ref{l:half_shannon},
where one uses Lemma~\ref{lem:disc-max2} in place of Lemma~\ref{lem:disc-max}.
We skip it.

For the proof of Theorem~\ref{thm:Bnorm-eq2}, one proceeds as in the proof
of Theorem~\ref{thm:Bnorm-eq}, using Lemma~\ref{lem:half-shannon2}
instead of Lemma~\ref{l:half_shannon}. The proof will be omitted.

\section{Application of weighted Besov spaces to nonlinear approximation}
\label{Nonlin-app}
\setcounter{equation}{0}

We consider here nonlinear n-term approximation for a needlet system
$\{\psi_\eta\}_{\eta\in \cX}$ with $\ph_\eta=\psi_\eta$,
defined as in (\ref{def.Phi-j})-(\ref{def-needlets2})
with $\hb = \ha$, $\ha\ge 0$. Then $\ha$ satisfies
$$
\ha^2(t)+\ha^2(2t)=1, \quad t\in [1/2, 1].
$$
Hence $\{\psi_\eta\}$ are real-valued.

Denote by $\Sigma_n$
the nonlinear set consisting of all functions $g$ of the form
$$
g = \sum_{\xi \in \Lambda} a_\xi \psi_\xi,
$$
where $\Lambda \subset \cX$, $\#\Lambda \le n$,
and $\Lambda$ is allowed to vary with $g$.
Let $\sigma_n(f)_p$ denote the error of best $\Lpp$-approximation to
$f \in \Lpp$ from $\Sigma_n$:
$$
 \sigma_n(f)_p := \inf_{g \in \Sigma_n} \|f - g\|_p.
$$
The approximation will take place in $\Lpp$, $0<p<\infty$.
Assume in the following that
$0 < p < \infty$, $s > 0$, and ${1}/{\tau} := s + {1}/{p}$.
Denote briefly $\BB s\tau:= \BBes s\tau\tau$.

By Theorem~\ref{thm:Bnorm-eq2} and (\ref{norm-Needlets})
it follows that
\begin{equation}\label{Btau-norm}
\|f\|_{\BB s\tau}\approx
\Big(\sum_{\xi\in\cX} \|\langle f, \psi_\xi \rangle \psi_\xi\|_p^\tau\Big)^{1/\tau}.
\end{equation}

The embedding of $\BB s\tau$ into $\Lpp$ plays an important role here.


\begin{proposition}\label{prop:embed-Lp}
If $f \in \BB s\tau$, then $f$ can be identified as a function
$f\in \Lpp$ and
\begin{equation}\label{embedding}
\|f\|_p \le
\Big\|\sum_{\xi\in\cX}|\langle f, \psi_\xi\rangle\psi_\xi(\cdot)|\Big\|_p
\le c\|f\|_{\BB s\tau}.
\end{equation}
\end{proposition}

We now state our main result in this section.


\begin{theorem}\label{t:jackson} {\bf [Jackson estimate]}
If $f \in \BB s\tau$, then
\begin{equation}\label{jackson}
\sigma_n(f)_p \le c n^{-s}\|f\|_{\BB s\tau},
\end{equation}
where $c$ depends only on $s$, $p$,
and the parameters of the needlet system.
\end{theorem}

The proofs of this theorem and Proposition~\ref{prop:embed-Lp}
can be carried out exactly as the proofs
of the respective Jackson estimate and embedding result
in \cite{NPW} and will be omitted.

It is an open problem to prove the companion to (\ref{jackson})
Bernstein estimate:
\begin{equation}\label{bernstein}
\|g\|_{\BB s\tau} \le c n^s \|g\|_p
\quad \hbox{for} ~~~g \in \Sigma_n,
\quad 1<p<\infty.
\end{equation}
This would enable one to characterize the rates (approximation spaces)
of nonlinear n-term approximation in $\Lpp$ ($1<p<\infty$) from needlet systems.

\section{Appendix}\label{appendix}
\setcounter{equation}{0}


\noindent
{\bf Proof of Proposition~\ref{thm:Lip}.}
We need the following integral representation of $\Ln(x, y)$
from \cite{PX1} (see (2.15)):
\begin{equation} \label{eq:Ln}
 \Ln(x, y) = c_{\a,\b} \int_0^\pi \int_0^1
    L_n^{\a,\b}(t(x,y,r,\psi)) \dm(r,\psi),
\end{equation}
where
$L_n^{\a,\b}(t)$ is defined by (\ref{def.Ln}), 
$$
 t(x,y,r,\psi) := \tfrac{1}{2}(1+x)(1+y)+ \tfrac{1}{2} (1-x)(1-y) r^2 +
      r\sqrt{1-x^2}\sqrt{1-y^2}  \cos \psi -1,
$$
the integral is against
$$
  \dm(r,\psi) := (1-r^2)^{\a-\b-1} r^{2\b+1} (\sin \psi)^{2\b} dr d\psi,
$$
and the constant $c_{\a,\b}$ is determined from
$$
c_{\a,\b}\int_0^\pi\int_0^1 1 \,\dm(r,\psi)=1.
$$

For any $u\in [-1, 1]$ we will denote by  $\theta_u$ the only angle in $[0, \pi]$
such that $u=\cos \theta_u$.

We will need the following lemma contained in the proof of
Theorem~2.4 in \cite{PX1}.


\begin{lemma}\label{lem1}
Let $\a,\b > -1/2$ and $k \ge 2\a+2\b+3$. Then there is a constant $c_k > 0$
depending only on $k$, $\a$, and $\b$ such that
for $x, y\in [-1, 1]$
$$
\int_0^\pi \int_0^1 \frac{n^{2 \alpha+1}\dm(r,\psi)}
{\left(1+ n \sqrt{1- t(x,y,r,\psi)}\right)^k}
\le c_k \frac{1} { {\sqrt{\W_{\a,\b}(n; x)}\sqrt{\W_{\a,\b}(n; y)}
(1+n|\t_x - \t_y|)^{\sigma}} },
$$
where $\sigma = k -2\a-2\b-3$.
\end{lemma}

Identity \eqref{eq:Ln} yields
\begin{equation}
\begin{aligned}\label{lip-1}
& |\Ln(x,y) - \Ln(\xi,y)| \cr
& \qquad  \le c \int_0^\pi \int_0^1
|L_n^{\a,\b}(t(x,y,r,\psi)) -L_n^{\a,\b}(t(\xi,y,r,\psi))|\dm(r,\psi)\cr
& \qquad  \le c \int_0^\pi \int_0^1
\| \partial L_n^{\a,\b}(\cdot) \|_{L^\infty(I_{r,\psi})}
|t(x,y,r,\psi) - t(\xi,y,r,\psi)| \dm(r,\psi),
\end{aligned}
\end{equation}
where $\partial f = f'$ and $I_{r,\psi}$ is the interval with end points
$t(x,y,r,\psi)$ and $t(\xi,y,r,\psi)$.

>From estimate (2.16) in \cite{PX1} and Markov's inequality, for any $k$
there exists a constant $c_k>0$ such that
\begin{align} \label{partialL}
& \|\partial L_n^{\a,\b}(\cdot)\|_{L^\infty(I_{r,\psi})}
\le c_k  \max_{u \in I_{r,\psi}}
\frac{n^{2\a + 4}}{\Big(1+n\sqrt{1-u}\Big)^k}\notag\\
& \qquad \le c_k  n^{2\a+4} \left[\Big(1+ n \sqrt{1-t(x,y,r,\psi)}\Big)^{-k} +
   \Big(1+ n \sqrt{1-t(\xi,y,r,\psi)}\Big)^{-k} \right].
\end{align}
For the rest of the proof we assume that $k>0$ is sufficiently large.

>From the definition of $t(x,y,r,\psi)$ one easily obtains
$$
 1-t(x,y,r,\psi)=  2 \sin^2 \frac{\t_x -\t_y}{2}
     + 2 \sin^2 \frac{\t_x}{2}\sin^2 \frac{\t_y}{2} (1-r^2) +
        \sin \t_x \sin\t_y (1- r \cos \psi),
$$
which implies
\begin{align*}
&t(x,y,r,\psi)-t(\xi,y,r,\psi) = \cos (\t_\xi-\t_y) - \cos(\t_x-\t_y)\\
& \qquad + (\cos\t_\xi - \cos\t_x) \sin^2 \frac{\t_y}2 (1-r^2) +
   (\sin \t_\xi - \sin\t_x) \sin \t_y (1-r \cos \psi).
\end{align*}
It is readily seen that
\begin{align*}
|\cos(\t_\xi-\t_y) - \cos(\t_x-\t_y)|
&= 2\Big|\sin\frac{\t_x+\t_\xi-2\t_y}{2} \sin \frac{\t_\xi-\t_x}{2}\Big| \\
&\le |\t_\xi - \t_x| (|\t_z - \t_y| + c n^{-1}),
\end{align*}
where we used that $|\t_z - \t_\xi|\le c n^{-1}$ and $|\t_x - \t_\xi|\le
c n^{-1}$.
Therefore,
\begin{align} \label{t-diff}
|t(x,y,r,\psi)&-t(\xi,y,r,\psi)|\notag\\
& \le |\t_\xi - \t_x|  \big[(|\t_z-\t_y| + c n^{-1})
+ \sin^2 \frac{\t_y}2 (1-r^2)
+  \sin \t_y (1-r \cos \psi) \big].
\end{align}
We use this and \eqref{partialL} in (\ref{lip-1}) to obtain
$$
|\Ln(x,y)-\Ln(\xi,y)| \le c|\t_\xi - \t_x|
        (A_1 + B_1 + A_2 + B_2 + A_3 +B_3),
$$
where $A_j$ and $B_j$ are integrals of the same type
with $A_j$ involving $t(x,y,r,\psi)$ and
$B_j$ involving $t(\xi,y,r,\psi)$;
the indices $j=1, 2, 3$ correspond to the three terms in the right-hand side
of \eqref{t-diff}.
We will estimate them separately.


\medskip\noindent
{\bf Case 1.} We first estimate the integral
$$
A_1: =  n^{2\a+4} \int_0^\pi \int_{0}^1 \frac{ |\t_z-\t_y| + c n^{-1} }
     { (1+ n \sqrt{1-t(x,y,r,\psi)})^{k}}\,\dm(r,\psi)
$$
and the integral $B_1$, the same as $A_1$ but
involving $t(\xi,y,r,\psi)$ in place of $t(x,y,r,\psi)$.

Using the estimate in Lemma \ref{lem1} and the fact that
$|\t_z - \t_y| \sim |\t_x - \t_y| + c n^{-1}$, we have
\begin{align*}
 A_1 & \le c \frac{ n^3 (|\t_z-\t_y| + c n^{-1})}
      { {\sqrt{\W(n; x)}\sqrt{\W(n; y)}
            (1+n|\t_x - \t_y|)^{\sigma}} } \\
 & \le c  \frac{ n^2 }
      { {\sqrt{\W(n; x)}\sqrt{\W(n; y)}
            (1+n|\t_z - \t_y|)^{\sigma -1}} }.
\end{align*}
On account of (\ref{omega<omega}) this gives the desired estimate.

The integral $B_1$ is estimated similarly with the same bound.


\medskip\noindent
{\bf Case 2.} We now estimate the integrals
$$
A_2: =  n^{2\a+4} \int_0^\pi \int_{0}^1 \frac{\sin^2 \frac{\t_y}2 (1-r^2)}
     { (1+ n \sqrt{1-t(x,y,r,\psi)})^{k}}\, \dm(r,\psi)
$$
and $B_2$ which is the same but involves $t(\xi,y,r,\psi)$
in place of $t(x,y,r,\psi)$.

By the definition of $dm_{\a,\b}(r,\psi)$ we have
$(1-r^2) dm_{\a,\b}(r,\psi) = dm_{\a+1,\b}(r,\psi)$.
Then using the estimate from Lemma \ref{lem1} with
$\alpha$ replaced by $\alpha +1$, we get
\begin{align*}
A_2 & \le c \frac{ n \sin^2 \frac{y}{2}}
{{\sqrt{\W_{\a+1,\b}(n; x)}\sqrt{\W_{\a+1,\b}(n; y)}
(1+n|\t_x - \t_y|)^{\sigma}} } \\
& \le c \frac{ n^2}
{{\sqrt{\W_{\a,\b}(n; x)}\sqrt{\W_{\a,\b}(n; y)}
(1+n|\t_x - \t_y|)^{\sigma}} },
\end{align*}
where we used the fact that
$\W_{\a+1,\b}(n; y) = \W_{\a,\b}(n; y)
(\sin^2 \frac{\t_y}2 + n^{-2})$ and hence
$\W_{\a+1,\b}(n; x) \ge \W_{\a,\b}(n; x) n^{-2}$.
The equivalence
$|\t_z - \t_y| \sim |\t_x - \t_y| + c n^{-1}$ and (\ref{omega<omega})
then give the desired estimate. The integral $B_2$ is estimated
similarly.


\medskip\noindent
{\bf Case 3.} We finally estimate the integrals
$$
A_3: =  n^{2\a+4} \int_0^\pi \int_{0}^1 \frac{\sin \t_y (1-r \cos \psi)}
     { (1+ n \sqrt{1-t(x,y,r,\psi)})^{k}} \, \dm(r,\psi)
$$
and $B_3$ which involves $t(\xi,y,r,\psi)$ in place of $t(x,y,r,\psi)$.

Assume first that $|\sin \t_x| \ge n^{-1}$. Using the fact that
$$
   1- t(x,y,r,\psi) \ge \sin \t_x \sin \t_y (1-r \cos\psi),
$$
we conclude that
\begin{align*}
A_3 &\le  \frac{n^{2\a+2}}{ \sin \t_x} \int_0^\pi \int_{0}^1 \frac{1}
     { (1+ n \sqrt{1-t(x,y,r,\psi)})^{k-2}} \, \dm(r,\psi)\\
    &\le c n^{2\a+3} \int_0^\pi \int_{0}^1 \frac{1}
     { (1+ n \sqrt{1-t(x,y,r,\psi)})^{k-2} } \,\dm(r,\psi).
\end{align*}
Now the estimate from Lemma~\ref{lem1} can be applied to get the desired
estimate.

Let $|\sin \t_x| \le n^{-1}$. We have
$$
 |\sin\t_y| \le | \sin\t_y - \sin\t_x| + |\sin \theta_x|
    \le |\t_y - \t_x| + n^{-1}
$$
and use the fact that
$$
 1- t (x,y,r,\psi) \ge 2 \sin^2 \frac{\t_x - \t_y}{2}\ge c(\t_x - \t_y)^2
$$
to conclude that
$$
  A_3 \le cn^{2\a+3} \int_0^\pi \int_{0}^1 \frac{1}
     { (1+ n \sqrt{1-t(x,y,r,\psi)})^{k-2}} \, d m_{\a,\b}(r,\psi).
$$
Applying the estimated from Lemma \ref{lem1} we obtain the desired result.
$B_3$ is estimated in the same way.

Putting the above estimates together completes the proof of Theorem~\ref{thm:Lip}.
$\qed$

\medskip


\noindent
{\bf Proof of Proposition~\ref{prop:lower-bound}.}
Note first that it suffices to prove (\ref{prop:lower-bound}) only for
$n\ge n_0$, where $n_0$ is sufficiently large. This follows from the fact
that $P_n^{(\a,\b)}$ and $P_{n+1}^{(\a,\b)}$ do not have common zeros and
$\omab(n; x) \sim 1$ if $n\le {\rm constant}$. Furthermore, since
$P_k^{(\a,\b)}(-x) = (-1)^k P_k^{(\b,\a)}(x)$, it is sufficient to consider
only the case $x \in [0,1]$.

Note that the Jacobi polynomials are normalized by
$P_k^{(\a,\b)}(1)= \binom{k+\a}{k}\sim k^\a$
and using Markov's inequality it follows that
$P_k^{(\a,\b)}(x)\ge ck^\a$ for $1-\delta k^{-2} \le x \le 1$,
where $\delta >0$ is a sufficiently small constant.
>From this one readily infers that (\ref{lowerbd}) holds
for $1-\delta_1 n^{-2} \le x \le 1$, $\delta_1>0$.
Define $\theta \in [0, \pi]$ from $x=\cos \theta$.
Then the latter condition on $x$ is apparently equivalent to
$0 \le \theta \le \delta_2 n^{-1}$ with $\delta_2$ being
a positive constant.

To estimate $\Lambda_n(\cos \theta)$ for $c^*n^{-1} \le \theta\le \pi/2$
with $c^*>0$ sufficiently large, we need the following asymptotic formula
of the Jacobi polynomials:
For $\a, \b>-1$,
\begin{align*}
 \left(\sin \frac{\t}{2}\right)^\a \left(\cos \frac{\t}{2}\right)^\b
 P_n^{(\a,\b)}(\cos\t)
 &= N^{-\a}\frac{\Gamma(n+\a+1)}{n!}\Big(\frac{\theta}{\sin \theta}\Big)^{1/2}
J_\alpha(N\t)\\
&\qquad+\theta^{1/2}\CO(n^{-3/2})
\end{align*}
if $c_0n^{-1}\le \theta\le \pi/2$,
where $N = n + \eta$ with $\eta:=(\a+\b+1)/2$, $J_\alpha$ is the Bessel
function, and $c_0>0$ is an arbitrary but fixed constant
(see \cite[Theorem 8.21.12, p. 195]{Sz}).

Using that
$2/\pi \le \sin \theta /\theta \le 1$ and
$\cos \theta/2 \sim 1$ on $[0, \pi/2]$,
and also $\Gamma(n+\a+1)/n!\sim n^\a$,
we infer from above
$$
\left(\sin \frac{\t}{2}\right)^{2\a} 
\left[ P_k^{(\a,\b)}(\cos\t) \right]^2 \ge
c_1[J_\alpha((k+\eta)\t)]^2
-c_2k^{-3/2}\theta^{1/2}|J_\alpha((k+\eta)\t)|.
$$
Recall the well-known asymptotic formula
$$
J_\alpha (z) = \left(\frac{2}{\pi z} \right)^{1/2}
\left[\cos (z + \gamma) + O(z^{-1})\right],
\quad z\to \infty,
$$
where $\gamma = -\alpha \pi/2 -\pi /4$.

All of the above leads to
\begin{align}\label{est-Lam-n}
\left(\sin \frac{\t}{2}\right)^{2\a}
\Lambda_n(\cos \theta)
&\ge \sum_{k= n}^{n+[\eps n]}
\Big( c_1[J_\alpha((k+\eta)\t)]^2\notag
-c_2k^{-3/2}\theta^{1/2}|J_\alpha((k+\eta)\t)|\Big)\\
&\ge  \frac{c}{n \t} \sum_{k= n}^{n+[\eps n]}
\left[\cos^2( k\t + b(\t)) - c'(n\theta)^{-1}\right]-c''\eps n^{-1},
\end{align}
for $c_0n^{-1}\le \theta\le \pi/2$,
where $b(\t) = ((\a+\b+1)/2)\t + \gamma$.
We now use the well known identities
for the Dirichlet kernel and its conjugate
to obtain, for $m>n$,
$$
\sum_{k=n}^m \cos^2(k\theta+b)
=\frac{1}{2}(m-n+1) +
(\cos 2b+\sin 2b) \frac{\sin(m-n+1)\theta\cos(n+m)\theta}{2\sin \theta}.
$$
Therefore,
$$
\sum_{k= n}^{n+[\eps n]}\cos^2 ( k \t + b(\t))
\ge \frac{1}{2}([\eps n]+1)\Big(1-\frac{2}{([\eps n]+1)\sin\theta}\Big)
\ge \frac{1}{2}\eps n\Big(1-\frac{\pi}{\eps n\theta}\Big)
\ge \frac{\eps n}{4},
$$
whenever $(2\pi/\eps)n^{-1}\le \theta\le \pi/2$.
Substituting this in (\ref{est-Lam-n}) we obtain
\begin{align} \label{est-Pk2}
\left(\sin \frac{\t}{2}\right)^{2\a}
\Lambda_n(\cos \theta)
&\ge  \frac{c}{n \t}\Big(\frac{\eps n}{4}- \frac{c'\eps n}{n\theta}\Big)
- \frac{c''\eps}{n} \\
& \ge \frac{c\eps}{\t}\left(\frac{1}{4} -\frac{c'}{c^*} \right)
   - \frac{c''\eps}{n}
\ge \frac{c^\diamond}{\theta}, \notag
\quad c^\diamond >0,
\end{align}
if $c^*n^{-1}\le \theta\le\pi/2$ with $c^*:=\max\{c_0, 8c', 2\pi/\eps\}$
and $n$ is sufficiently large.
Hence,
\begin{align}\label{Lambda1}
\Lambda_n(\cos \theta)
\ge \frac{c^\diamond}{\theta} \left(\sin \frac{\t}{2}\right)^{-2\a}
\ge c\theta^{-(2\a+1)},
\quad c^*n^{-1}\le \theta\le\pi/2,
\end{align}
for sufficiently large $n$,
which yields (\ref{lowerbd}) in this case.

\medskip

For the remaining case $\delta_2n^{-1}\le \theta \le c^* n^{-1}$, we need
further properties of Jacobi polynomials.
Let $x_{\nu,n} = \cos \t_{\nu,n}$ denote the
zeros of Jacobi polynomial $P_n^{(\a,\b)}$, where
$$
   0 < \t_{1,n} < \t_{2,n} < \cdots < \t_{n,n}  < \pi.
$$
It is well known that
$\t_{\nu,n}\sim \nu/n$,
but we will need much more precise asymptotic representation for $\t_{\nu,n}$,
see below.
The Jacobi polynomials satisfy the following relation
(see e.g. \cite[Theorem 3.3, p. 171]{N}),
\begin{equation} \label{JacobiEst}
P_n^{(\a,\b)}(\cos \theta) \sim n^{1/2} |\t - \t_{\nu_\t,n}|
\frac{n^{\a+1/2}}{\nu_\t^{\a+1/2}},
\quad \t \in [0,\pi],
\end{equation}
where $\nu_\t$ denotes the index of the zero $x_{\nu,n}$,
$1\le \nu \le n$,
which is (one of) the closest to $x$ ($x=\cos \theta$).

We will need the asymptotics of the zeros of the Jacobi polynomials from
\cite{FW}:
\begin{equation} \label{asymp-zeros}
\t_{\nu,n} = \frac{j_{\a,\nu}}{N} + \frac{1}{4 N^2}
\left[\left(\a^2-\frac{1}{4}\right)\frac{1-t\cot t}{2t}
-\frac{\a^2-\b^2}{4}\tan\frac{t}{2} \right]+ t^2 \CO(n^{-3}),
\end{equation}
where $N=n + \eta$ as before,
$j_{\alpha,\nu}$ is the $\nu$th positive zero of
the Bessel function $J_\alpha(x)$
and $t = j_{\alpha,\nu} /N$.
Here the $\CO$-term is uniformly bounded for $\nu =1,2,\ldots, [\gamma n]$,
where $\gamma\in (0, 1)$.
It is easy to verify that
$(1-t\cot t)/t = \CO(t)$ as $t \to 0$
and obviously $1/(n+\eta)- 1/n = \CO(n^{-2})$.
Hence
\begin{equation}\label{zero}
\theta_{\nu,n} = \frac{j_{\a,\nu}}{n} + \CO(n^{-2}),
\quad \nu =1,2,\ldots, [\gamma n].
\end{equation}

We will also use that
$$
0<j_{\alpha, 1}< j_{\alpha, 2}<\cdots
\quad \mbox{ and } \quad
j_{\alpha,\nu} \to \infty.
$$
Let $j_{\alpha, \nu_{\max}}:=\max\{j_{\alpha,\nu}:
  j_{\alpha,\nu}\le (1+\eps)c^*\}$
and denote
$\CJ:= \{ j_{\alpha, 1}, j_{\alpha, 2}, \dots, j_{\alpha, \nu_{\max}}\}$.
Notice that $\nu_{\max}$ is a constant independent of $n$.
Suppose that $\CJ\ne\emptyset$ (the case $\CJ=\emptyset$ is easier).

Fix $\delta_2n^{-1}\le \theta \le c^* n^{-1}$.
Then by \eqref{JacobiEst} it follows that
$$
P_k^{(\a,\b)}(\cos \theta) \sim k^{\a +1} |\t - \t_{\nu_\t,k}|,
$$
where the $\nu_\theta$'s involved are bounded by a constant independent of
$n$. Hence, (\ref{zero}) can be used to represent $\t_{\nu_\t,k}$ for
$n\le k\le n+[\eps n]$ if $n$ is sufficiently large. Using the above we get
\begin{align*}
\Lambda_n(\cos \theta)
&\ge c n^{2 \a+2} \sum_{k=n}^{n+[\eps n]}|\t - \t_{\nu_\t,k}|^2
\ge c n^{2 \a} \sum_{k=n}^{n+[\eps n]}|k\t - k\t_{\nu_\t,k}|^2\\
& \ge c n^{2 \a} \sum_{k=n}^{n+[\eps n]}
\Big(|k\t - j_{\a, \nu_\t}|^2-c'k^{-1}|k\t - j_{\a, \nu_\t}|\Big)\\
& \ge c n^{2 \a} \Big(\sum_{k=n}^{n+[\eps n]}
|k\t - j_{\a, \nu_\t}|^2-c'c^*\eps\Big),
\end{align*}
where we used (\ref{zero}).
Therefore,
\begin{align}\label{sum1}
\Lambda_n(\cos \theta)
& \ge c n^{2 \a}\Big(\sum_{k=n}^{n+[\eps n]}\dist (k\t, \CJ)^2
  - c^\diamond\Big),
\quad c, c^\diamond >0,
\end{align}
where $\dist (k\t, \CJ)$ denotes the distance of $k\t$ from the set
$\CJ$, that is the distance of $k\t$ from the nearest zero of the Bessel
function $J_\a(x)$.

It remains to estimate the sum in (\ref{sum1}).
Denote $\CK:=\{n, n+1, \dots, n+[\eps n]\}$ and
let $\CK_0$ be the set of all indices $k\in \CK$ such that
$\dist (k \theta, \CJ) <m\theta$,
where $m:=[\eps n/(6\nu_{\max})]$.
Evidently
$$
\# \CK_0 \le (2m+1)\nu_{\max} \le (2[\eps n/(6\nu_{\max})]+1)\nu_{\max}
\le \eps n/2
\quad\mbox{if $n \ge 6\nu_{\max}\eps^{-1}$.}
$$
Then
$\# \CK \setminus \CK_0 \ge [\eps n]+1 - \eps n/2 \ge \eps n/2$
and hence
$$
\sum_{k=n}^{n+[\eps n]}\dist (k\t, \CJ)^2
\ge \sum_{k\in\CK\setminus \CK_0}(m\theta)^2
\ge c\sum_{k\in\CK\setminus \CK_0}(n\theta)^2
\ge c\delta_2^2\eps n\ge c_*n,
\quad c_*>0.
$$
Inserting this in (\ref{sum1}) we obtain
$$
\Lambda_n(\cos \theta)
\ge c n^{2 \a}(c_*n - c^\diamond)
\ge \tilde c n^{2 \a+1}
$$
for sufficiently large $n$.
This implies the stated inequality \eqref{lowerbd} with $x=\cos \theta$
in the case $\delta_2n^{-1}\le \theta \le c^* n^{-1}$.
The proof of Proposition~\ref{prop:lower-bound} is complete.
$\qed$

\medskip


\medskip\noindent
{\bf Proof of Proposition~\ref{Nikolski}.}
Suppose $\alpha \ge \beta$ and let $1<q<\infty$.
By Lemma~\ref{lem:Ker-n}, (i) we have $g=L_n*g$ and using
H\"older's ineqaulity, (\ref{est-Lp-int}), and that
$\W(n;x)\ge cn^{-2\a-1}$ we obtain
$$
|g(x)| \le \|g\|_q\left(\frac{n}{\W(n;x)}\right)^{1/q}
\le c n^{(2 \a +2)/q}\|g\|_q,
\quad x\in [-1, 1],
$$
which leads to
\begin{equation}\label{Nik1}
\|g\|_\infty \le c n^{(2\a+2)/q}\|g\|_q,
\quad 1<q\le \infty.
\end{equation}
If $0<q\le 1$, then the above inequality with $q=2$ gives
$$
\|g\|_\infty^2
\le c n^{2\a+2} \int_{-1}^1 |g(y)|^{2-q} |g(y)|^q \w(y)dy
\le c n^{2\a+2}\|g\|_\infty^{2-q}\|g\|_q^q
$$
which shows that (\ref{Nik1}) holds for $0<q\le 1$ as well.

Let $0<q<p<\infty$ (the case $p=\infty$ is contained in (\ref{Nik1})).
Then using (\ref{Nik1}) we obtain
\begin{align*}
\|g\|_p  & = \left(\int_{-1}^1 |g(x)|^{p-q} |g(x)|^q \w(x) dx \right)^{1/p}\\
   & \le c n^{(2\a +2)(\frac{1}{q} - \frac{1}{p})} \|g\|_q^{\frac{p-q}{p} }
     \|g\|_q^{\frac{q}{p}}
    =  c n^{(2\a +2)(\frac{1}{q} - \frac{1}{p})} \|g\|_q.
\end{align*}
On the other hand, by \cite[p. 114]{N}
$$
 \|g\|_p \le c n^{2 (\frac{1}{q}-\frac{1}{p} )} \|g\|_q.
$$
Putting the above two estimates together gives \eqref{norm-relation}.


To prove (\ref{norm-relation2}) we will need the following inequality
\begin{equation} \label{aaa}
 \int_{-1}^1 \frac{\w(y)}{\W(n;y)^{p/2 + \gamma} (1 + n d(x,y))^\sigma}dy
    \le c \frac{1}{n \W(n;x)^{p/2 + \gamma -1}},
\quad x\in [-1, 1],
\end{equation}
where $\gamma \in \RR$ and $\sigma$ is sufficiently large.
The proof of \eqref{aaa} is contained in the proof of
Proposition 1 in \cite{PX1}.
Assume $1<q<\infty$. Then using Lemma~\ref{lem:Ker-n}, (i),
H\"older's inequality ($1/q + 1/q'=1$), and (\ref{Lbound1})
we have, for $x\in[-1, 1]$,
\begin{align*}
|g(x)| & \le \|\W(n;\cdot)^{s +\frac{1}{p} - \frac{1}{q}} g(\cdot)\|_q
\left( \int_{-1}^1 \Big|L_n(x,y) \W(n;y)^{-s-\frac{1}{p}+\frac{1}{q}}\Big|^{q'}
\w(y) dy \right)^{1/q'} \\
& \le c \frac{n} {\W(n;x)^{1/2}} \left (\int_{-1}^1  \frac{\w(y)dy}
{ \W(n;y)^{\frac{q'}{2} + \gamma} (1+ n d(x,y))^\sigma } \right)^{1/q'}
\|\W(n;\cdot)^{s +\frac{1}{p} - \frac{1}{q}} g(\cdot)\|_q
\end{align*}
with $\gamma = q'(s + \frac{1}{p}-\frac{1}{q})$.
Now, applying \eqref{aaa} we infer
$$
 |g(x)|  \le  c \frac{n^{1/q}}{\W(n;x)^{s +1/p}}
     \|\W(n;\cdot)^{s +\frac{1}{p} - \frac{1}{q}} g(\cdot)\|_q,
$$
which implies
\begin{equation} \label{aaa2}
\|\W(n;\cdot)^{s+1/p}g(\cdot)\|_\infty \le c n^{1/q}
\|\W(n;\cdot)^{s +\frac{1}{p} - \frac{1}{q}} g(\cdot)\|_q,
\quad 1<q\le \infty.
\end{equation}
If $0 < q \le 1$, then by (\ref{aaa2}) with $q=2$
\begin{align*}
 \|\W(n;\cdot)^{s+1/p}g(\cdot)\|_\infty  &\le c n^{1/2}
 \|\W(n;\cdot)^{s +\frac{1}{p} - \frac{1}{2}} g(\cdot)\|_2 \\
 & \le c n^{1/2}  \|\W(n;\cdot)^{s+1/p}g(\cdot)\|_\infty^{1-q/2}
 \|\W(n;\cdot)^{s +\frac{1}{p} - \frac{1}{q}} g(\cdot)\|_q^{q/2},
\end{align*}
which shows that \eqref{aaa2} holds for $0 < q \le 1$ as well.
Let $p<\infty$.
Using (\ref{aaa2}), we get
\begin{align*}
\|\W(n;\cdot)^s g(\cdot)\|_p & = \left( \int_{-1}^1 \left[\W(n;x)^s g(x)\right]^{p-q}
\left[\W(n;x)^s g(x)\right]^{q} \w(x) dx \right)^{1/p} \\
\le & c n^{\frac{1}{q} - \frac{1}{p}}
\|\W(n;\cdot)^{s +\frac{1}{p} - \frac{1}{q}} g(\cdot)\|_q^{p-q}
\left( \int_{-1}^1 \frac{\left[\W(n;x)^s g(x)\right]^{q}}
{\W(n;x)^{\frac{p-q}{p}}}  \w(x) dx \right)^{1/p} \\
= & c n^{\frac{1}{q}-\frac{1}{p}}
\|\W(n;\cdot)^{s +\frac{1}{p} - \frac{1}{q}} g(\cdot)\|_q.
\end{align*}
Thus \eqref{norm-relation2} is established if $p < \infty$.
When $p = \infty$ \eqref{norm-relation2} follows from \eqref{aaa2}.
$\qed$

\medskip


\noindent
{\bf Proof of estimate (\ref{omega<omega}).}
We only consider the case when $-1/2 \le x, y\le 1$ since the other cases
are similar or simpler.
Choose $\theta, \phi \in [0, 2\pi/3]$ so that $x=\cos \theta$, $y=\cos \phi$.
Then $d(x, y)=|\theta-\phi|$.
We have
\begin{align*}
\sin\theta+n^{-1}
&\le |\sin\theta-\sin\phi|+\sin\phi+n^{-1}
\le |\theta-\phi|+\sin\phi+n^{-1}\\
&\le (1+n|\theta-\phi|)(\sin\phi+n^{-1}).
\end{align*}
Then (\ref{omega<omega}) follows using (\ref{wn-sim1}).
$\qed$

\medskip


\noindent
{\bf Proof of Lemma~\ref{lem:J-maximal}.}
We first show that for $y\in[0,1 ]$ and $0<r\le\pi$
\begin{equation}\label{mu-B}
\mu(B_y(r)):= \int_{B_y(r)}\w(x)dx
\sim r(d(y, 1)+r)^{2\a+1}.
\end{equation}
Indeed, choose $0\le \zeta\le \pi/2$ so that $y=\cos\zeta$
and consider the case when $\zeta+r\le 2\pi/3$
(the case $\zeta+r>2\pi/3$ is trivial; then $\mu(B_y(r))\sim 1$).
We have
\begin{align*}
\mu(B_y(r))
&\sim \int_{\max\{\zeta-r, 0\}}^{\zeta+r}(1-\cos u)^\a\sin u du
\sim \int_{\max\{\zeta-r, 0\}}^{\zeta+r}u^{2\a+1}du\\
&\sim (\zeta +r-\max\{\zeta-r, 0\})(\zeta+r)^{2\a+1}
\sim r(\zeta+r)^{2\a+1},
\end{align*}
which yields (\ref{mu-B}).


We now proceed with the proof of (\ref{J-max1}).
Denote briefly $J_\eta:=B_\eta(\eps)$ which is an interval.
Equivalence (\ref{J-max1}) is trivial when $x\in J_\eta$.

Assume $x\in [-1, 1]\setminus J_\eta$.
Denote by $I_{x,\eta}$
the interval with end points $x$ and $\eta$.
By the definition of the maximal operator in (\ref{def.max-fun})
it readily follows that
$$
\left(\frac{\mu(J_\eta)}{\mu(J_\eta\cup I_{x,\eta})}\right)^{1/t}
\le(\cM_t\ONE_{J_\eta})(x)\le
\left(\frac{\mu(J_\eta)}{\mu(I_{x,\eta})}\right)^{1/t}
$$
and since $\mu(I_{x,\eta})\le \mu(J_\eta\cup I_{x,\eta})
\le \mu(J_\eta)+\mu( I_{x,\eta})\le c \mu(I_{x,\eta})$,
we have
\begin{equation}\label{sim-maximal}
 (\cM_t\ONE_{J_\eta})(x) \sim \left(\frac{\mu(J_\eta)}{\mu(I_{x,\eta})}\right)^{1/t}.
\end{equation}

We will only consider the case when $x\in[-1/2, 1]$,
since the case $x\in[-1, -1/2]$ is simpler.
Choose $\gamma\in [0, \pi/2]$ and $\phi\in [0, 2\pi/3]$
so that $\eta=\cos \gamma$ and $x=\cos \phi$.
Then $d(\eta, x)=|\gamma-\phi|$.
By (\ref{mu-B}) and (\ref{sim-maximal}) it follows that
$$
[(\cM_t\ONE_{J_\eta})(x)]^t
\sim \frac{\eps(d(\eta, 1)+\eps)^{2\a+1}}
{\frac{|\gamma-\phi|}{2}\Big(\frac{\gamma+\phi}{2}+\frac{|\gamma-\phi|}{2}\Big)^{2\a+1}}
\sim \frac{\eps}{d(\eta, x)}
\Big(\frac{d(\eta, 1)+\eps}{d(\eta, x)+d(\eta, 1)}\Big)^{2\a+1},
$$
which implies (\ref{J-max1}).
Estimates (\ref{J-max3}) are immediate from (\ref{J-max1}).
$\qed$

\medskip


\noindent
{\bf Proof of (\ref{norm-Needlets})-(\ref{norm-Needlets2}).}
The equivalence
$\norm{\tONE_{I_\xi}}_\Lp\sim (2^{-j}\W(2^j;\xi))^{1/p-1/2}$
follows by (\ref{muI}).

>From (\ref{est-Lp-int}) and (\ref{coeff}) it follows that,
for $\xi\in\cX_j$ and $0<p<\infty$,
\begin{equation}\label{upper-est}
\|\ph_\xi\|_\Lp, \|\psi_\xi\|_\Lp
\le c \cc_\xi^{1/2}\Big(\frac{2^j}{\W(2^j;\xi)}\Big)^{1-1/p}
\le c \Big(\frac{2^j}{\W(2^j;\xi)}\Big)^{1/2-1/p}.
\end{equation}
When $p=\infty$, similar estimates follow by (\ref{local-needlets3}).

To estimate $\|\ph_\xi\|_\Lp$, $\|\psi_\xi\|_\Lp $ from below,
we first note that by (\ref{est-Lp-norm}) and (\ref{upper-est}) it follows that
$\|\ph_\xi\|_\L2\sim\|\psi_\xi\|_\L2\sim 1$.
Let $2<p <\infty$ and $\frac 1p+\frac1{p'}=1$.  Using H\"older's inequality
and (\ref{upper-est}) we obtain
$$
0<c\le \|\ph_\xi\|_\L2^2
\le\|\ph_\xi\|_\Lp \|\ph_\xi\|_{p'}
\le\|\ph_\xi\|_\Lp\Big(\frac{2^j}{\W(2^j;\xi)}\Big)^{1/2-1/p'}
$$
and similarly for $\psi$.
Hence
\begin{equation}\label{lower-est}
\|\ph_\xi\|_\Lp, \|\psi_\xi\|_\Lp
\ge c \Big(\frac{2^j}{\W(2^j;\xi)}\Big)^{1/2-1/p}.
\end{equation}
In the case $p=\infty$, we proceed similarly as above and obtain the same estimate.

If $0<p<2$, then
\begin{align*}
0<c\le \|\ph_\xi\|_\L2^2
\le \|\ph_\xi\|_\Lp^p \|\ph_\xi\|_\Linfty^{2-p}
\le c\|\ph_\xi\|_\Lp^p \Big(\frac{2^j}{\W(2^j;\xi)}\Big)^{1-p/2},
\end{align*}
which implies (\ref{lower-est}) and similarly for $\|\psi_\xi\|_\Lp$.

Finally, (\ref{norm-Needlets2}) follows by the lower bound from (\ref{lower-est})
with $p=\infty$ and (\ref{local-Needlets21}).
$\qed$

\medskip

\noindent
{\bf Proof of Lemma~\ref{lem:Max-needl}.}
Estimate (\ref{Max-needl1}) is immediate from (\ref{local-needlets3})
and Lemma~\ref{lem:J-maximal} (see (\ref{J-max3}) and (\ref{muI})).

For the proof of (\ref{Max-needl2}) we first observe that by (\ref{norm-Needlets2})
there exists a point $\zeta\in B_\xi(c^*2^{-j})$ such that
$|\ph_\xi(\zeta)| \ge c^\diamond (2^j/\W(2^j;\xi))^{1/2}$.
By (\ref{Lip}) it follows that
$$
|\ph_\xi(x)- \ph_\xi(\zeta)| \le c\frac{2^{3j/2}d(x, \zeta)}{\sqrt{\W(2^j;\xi)}},
\quad x\in B_\zeta(c2^{-j}),
$$
and hence for sufficiently small constant $c^\flat>0$ we have
$$
|\ph_\xi(x)| \ge \frac{c^\diamond}{2} (2^j/\W(2^j;\xi))^{1/2}
\quad\mbox{for $x\in B_\zeta(c^\flat2^{-j})\cap[-1, 1]$.}
$$
Therefore, there exists an interval $J_\eta:= B_\eta(c2^{-j}) \subset [-1, 1]$
such that
$$
\|\ph_\xi\|_{L^\infty(J_\eta)} \ge c(2^j/\W(2^j;\xi))^{1/2}
\quad\mbox{and}\quad
d(\eta, \xi)\le c2^{-j}.
$$
Hence,
$$
(\cM_t\ph_\xi)(x) \ge c(\cM_t\tONE_{J_\eta})(x) \ge c\tONE_{I_\xi}(x),
$$
where for the last estimate we used (\ref{J-max1}).
Thus (\ref{Max-needl2}) is established for $\ph_\xi$.
The proof for $\psi_\xi$ is the same.
$\qed$

\medskip

\noindent
{\bf Proof of Lemma \ref{Almdiag}.}
>From the orthogonality of Jacobi polynomials it follows that
$\Phi_j\ast \psi_\xi(x)=0$ if $\xi\in\cX_\nu$, where
$\nu\ge j+2$ or $\nu\le j-2$.

Assume that $\xi \in \cX_\nu$, $j-1\le \nu\le j+1$.
>From the localization of the kernels $\Phi_j, \Psi_\nu$
(see (\ref{local-Needlets2})) and the definition of $\cc_\xi$
(see (\ref{quadrature1}))
we get that for any $\sigma>0$ there is a constant $c_\sigma>0$ such that
\begin{equation*}
\begin{aligned}
&|\Phi_j\ast \psi_\xi(x)|
=\sqrt{\cc_\xi}\Big|\int_{-1}^1\Phi_j(x,y)\Psi_\nu (y,\xi)\w(y)\, dy\Big|\\
&\le  c 2^{3j/2}\W(2^j; x)^{-1/2}\int_{-1}^1\frac{ \w(y)}
{ \W(2^j; y)(1+2^jd(x,y))^\sigma(1+2^jd(y,\xi))^\sigma}\, dy.
\end{aligned}
\end{equation*}
Setting $\xi=\cos \theta$,  $x=\cos \eta$ for some $0\le \theta,\eta\le \pi$
and applying the substitution $y=\cos \phi$,
we obtain
\begin{align*}
|\Phi_j\ast \psi_\xi(x)|
&\le  c2^{j3/2}\W(2^j; x)^{-1/2}\int_{0}^\pi
 \frac{\w(\cos \phi)\sin \phi}{\W(2^j; \cos \phi)
(1+2^j|\eta-\phi|)^\sigma(1+2^j|\theta-\phi|)^\sigma}\, d\phi\\
&\le c2^{j3/2}\W(2^j; x)^{-1/2}
\int_0^\pi\frac{1}{(1+2^j|\eta-\phi|)^\sigma(1+2^j|\theta-\phi|)^\sigma}\,d\phi\\
&\le  c2^{j/2}\W(2^j; x)^{-1/2}(1+2^j|\eta-\theta|)^{-\sigma},
\end{align*}
where we used the inequality
$$
2^j\int_\R \frac1{
(1+2^j|\eta-\phi|)^\sigma(1+2^j|\theta-\phi|)^\sigma}\, d\phi
\le c (1+2^j|\eta-\theta|)^{-\sigma}. \qed
$$

\medskip

For the proof of Lemma~\ref{l:weak_inequality} we will need this lemma.


\begin{lemma}\label{lem:quasi-Lip}
Let $P\in \Pi_{2^j},j \ge 0$ and $\xi\in \cX_j$.
Suppose $x_1,x_2\in [-1,1]$ and
$d(x_\nu,\xi)\le c_{\star}2^{-j}$, $\nu=1,2$.
Then for any $\sigma>0$
$$
|P(x_1)-P(x_2)|\le c_\sigma 2^j d(x_1, x_2)
\sum_{\eta\in\cX_j}\frac{|P(\eta)|}{(1+2^j d(\xi, \eta))^\sigma},
$$
where $c_\sigma>0$ depends only on $\sigma$, $\a$, $\b$, and $c_\star$.
\end{lemma}

\noindent
{\bf Proof.}
Let $P\in\Pi_{2^j}$.
Suppose $L_{2^j}(x,y)$ is the reproducing kernel from Lemma~\ref{lem:Ker-n}
with $n=2^j$. Then $L_{2^j}*P=P$.
Since $L_{2^j}(x, \cdot)P(\cdot) \in \Pi_{2^{j+2}}$, using that the quadrature
(\ref{quadrature1}) is exact for all polynomials from $\Pi_{2^{j+2}}$ we get
$$
P(x)=\int_{-1}^1L_{2^j}(x,y)P(y)\w(y)dy
=\sum_{\eta\in\cX_j}\cc_\eta L_{2^j}(x,\eta)P(\eta),
\quad x\in[-1, 1].
$$
Recall that $\cc_\eta\sim 2^{-j}\W(2^j;\eta)$.
Now, using Theorem~\ref{thm:Lip} we obtain
for $x_1, x_2\in[-1, 1]$ with $d(x_\nu, \xi)\le c_\star2^{-j}$, $\nu=1, 2$,
\begin{align*}
|P(x_1)-P(x_2)|
&=\Big|\int_{-1}^1[\Ker_{2^j}(x_1,y)-\Ker_{2^j}(x_2,y)]P(y)\w(y) \, dy\Big|\\
&\le\sum_{\eta\in \cX_j}
|\cc_{\eta}||\Ker_{2^j}(x_1,\eta)- \Ker_{2^j}(x_1,\eta)||P(\eta)|\\
&\le c2^j d(x_1,x_2) \sum_{\eta\in \cX_j}
\Big(\frac{\W(2^j;\eta)}{\W(2^j;\xi)}\Big)^{1/2}
\frac{|P(\eta)|}{(1+2^jd(\xi,\eta))^{\sigma}}\\
&\le c2^j d(x_1,x_2) \sum_{\eta\in \cX_j}
\frac{|P(\eta)|}{(1+2^jd(\xi,\eta))^{\sigma-\max\{\a,\b\}-1/2}},
\end{align*}
where for the last inequality we used (\ref{omega<omega}).
Since $\sigma>0$ can be arbitrarily large the result follows.
$\qed$

\medskip

\noindent
{\bf Proof of Lemma \ref{l:weak_inequality}.}
Clearly
$a_\xi\le b_\xi +d_\xi$,
where
$$
d_\xi:=\max\{|P(x_1)-P(x_2)|: x_1\in I_\xi, d(x_1,x_2)\le c_22^{-(j+r)}\},
$$
and $c_2$ is the constant appearing in (\ref{I-B}).
%
By Lemma \ref{lem:quasi-Lip} it follows that
$$
d_\xi\le c2^{-r} \sum_{\eta\in\cX_j}\frac{|P(\eta)|}{(1+2^j d(\xi, \eta))^\sigma},
\quad \xi\in \cX_j.
$$
Then recalling the definition of $d_\xi^*$ in (\ref{def.h-star}) we infer
\begin{align*}
d_\xi^\ast&\le c 2^{-r} \sum_{w\in \cX_j}\sum_{\eta\in\cX_j}
\frac{|P(\eta)|}{(1+2^j d(w,\eta))^\sigma (1+2^jd(\xi,w))^\sigma}\\
&\le c 2^{-r}\sum_{\eta\in \cX_j}\frac{|P(\eta)|}{(1+2^j d(\eta, \xi))^\sigma}
\le c 2^{-r}a_\xi^*,
\end{align*}
where for the second inequality we switched the order of summation and
used the simple fact that for $\sigma>1$
$$
\sum_{w\in \cX_j} \frac1
{  (1+2^j d(w,\eta))^\sigma (1+2^jd(\xi,w))^\sigma}
\le c\frac1{  (1+2^j d(\eta, \xi))^\sigma}.
$$
Consequently,
$a_\xi^*\le b_\xi^*+d_\xi^* \le b_\xi^*+c2^{-r}a_\xi^*$ with $c>0$
independent of $r$.
Choosing $r$ sufficiently large we obtain
$a_\xi^*\le cb_\xi^\ast$.
The estimate in the other direction is trivial.
$\qed$

\medskip


\noindent
{\bf Proof of Lemma \ref{lem:disc-max}.}
We may assume that $\a\ge \b$.
Fix $\xi\in\cX_j$ and define
$Y_0:=\{\eta\in \cX_j: d(\eta,\xi)\le c_22^{-j}\}$ and
$$
Y_m :=\{\eta\in \cX_j: c_22^{-j+m-1} < d(\eta,\xi)\le c_2 2^{-j+m}\}, \quad m\ge 1.
$$
where $c_2>0$ is from (\ref{I-B}).
Using (\ref{I-B}) we have $\# Y_m\le c2^m$.
Also, let
$$
J_m := B_\xi(c_2(2^m+1)2^{-j})
= \{x\in[-1,1]: d(\eta,\xi)\le c_2(2^m+1)2^{-j}\}, \quad m\ge 0.
$$
Evidently, $J_m$ is an interval and $I_\eta\subset J_m$ if $\eta\in Y_\nu$,
$0\le \nu\le m$.

We next show that
\begin{equation}\label{disc-max2}
\mu(J_m) \le c 2^{m(4\a +3)}\mu(I_\eta)
\quad\mbox{for all}\quad \eta\in Y_m.
\end{equation}
Suppose $\xi\in[0,1]$; the case $\xi\in[-1,0]$ is the same.
Let $J_m=:[y_1, y_2]$ and chose $\phi_1, \phi_2\in[0, \pi]$ so that
$y_1=:\cos \phi_1$ and $y_2=:\cos \phi_2$ ($\phi_1>\phi_2$).
Exactly as in the proof of Lemma~\ref{lem:J-maximal}
$$
\mu(J_m) = \int_{J_m}\w(y)dy\le c(\phi_1-\phi_2)\phi_1^{2\a+1}
\le c 2^{-j+m}\W(2^j, y_1)
$$
and using (\ref{omega<omega})
\begin{equation}\label{est-muJ}
\mu(J_m) \le c2^{-j+m}\W(2^j,\xi)(1+2^jd(y_1,\xi))^{2\a+1}
\le c2^{-j+(2\a+2)m}\W(2^j,\xi).
\end{equation}
On the other hand, using again (\ref{omega<omega})
$$
\mu(I_\eta) \ge c2^{-j}\W(2^j,\eta) \ge c 2^{-j}\W(2^j,\xi)(1+2^jd(\eta,\xi))^{-2\a-1}
\ge c 2^{-j-(2\a+1)m}\W(2^j,\xi).
$$
Combining this with (\ref{est-muJ}) gives (\ref{disc-max2}).

Let $\rho:=\max \{0, 1-\frac1{t}\}<1$.
Using H\"older's inequality if  $t> 1$ and
the $t$-triangle inequality if $0<t\le1$, we have
\begin{align*}
b_\xi^{\ast}=
\sum_{\eta\in\cX_j}\frac{|b_\eta|}{(1+2^jd(\eta,\xi))^\sigma}
\le c \sum_{m\ge 0}2^{-m\sigma}\sum_{\eta\in Y_m} |b_\eta|
\le c\sum_{m\ge 0}2^{-m(\sigma-\rho)}(\sum_{\eta\in Y_m} |b_\eta|^t)^{1/t}.
\end{align*}
We next use (\ref{disc-max2}) to obtain, for $x\in I_\xi$,
\begin{align*}
b_\xi^{\ast}
&=c \sum_{m\ge 0}2^{-m(\sigma-1)}\Bigl(\int_{-1}^1
\Big[\sum_{\eta\in Y_m}
|b_\eta|\mu(I_\eta)^{-1/t}\ONE_{I_\eta}(x)\Big]^t\w(x)\, dx\Bigr)^{1/t}\\
&\le c \sum_{m\ge 0}2^{-m(\sigma-1)}\Bigl(\frac1{\mu(J_m)}\int_{J_m}
\Big[\sum_{\eta\in Y_m}
\Big(\frac{\mu(J_m)}{\mu(I_\eta)}\Big)^{1/t}|b_\eta|\ONE_{I_\eta}(x)\Big]^t
\w(x)\,dx\Bigr)^{1/t}\\
&\le c\sum_{m\ge 0}2^{-m(\sigma-1-(4\a +3)/t)}
\Bigl(\frac1{\mu(J_m)}\int_{J_m}
\Big[\sum_{\eta\in Y_m}|b_\eta|\ONE_{I_\eta}(x)\Big]^t\w(x)\, dx\Bigr)^{1/t}\\
&\le c\cM_t\Big(\sum_{w\in \cX_j}|b_w|\ONE_{I_\omega}\Big)(x)
\sum_{m\ge 0}2^{-m(\sigma-1-(4\a +3)/t)}
\le c\cM_t\Big(\sum_{w\in \cX_j}|b_w|\ONE_{I_\omega}\Big)(x),
\end{align*}
where for the last inequality we used that
$\sigma>(4\a+3)/t+1$.
\qed

\medskip


\noindent
{\bf Proof of Lemma \ref{l:half_shannon}.}
For $\xi\in\cX_j$, we set
$a_\xi:=\max_{x\in I_\xi}|P(x)|$, $m_\xi:=\min_{x\in I_\xi}|P(x)|$ and
$$
b_\xi:=\max\{\min_{x\in I_w }|P(x)|:w\in \cX_{j+r},I_w\cap I_\xi\ne\emptyset \},
$$
where $r\ge 1$ (sufficient large) is the constant from Lemma \ref{l:weak_inequality}.
If $0<t<p$ then
\begin{equation}
\begin{aligned}\label{half_shannon}
\Big(\sum_{\xi\in \cX_j}a_\xi^p \mu(I_\xi)\Big)^{1/p}
&= \Big\|\sum_{\xi\in \cX_j}a_\xi \ONE_{I_\xi}(\cdot)\Big\|_{\Lp}
\le c\Big\|\sum_{\xi\in \cX_j}b_\xi^\ast \ONE_{I_\xi}(\cdot)\Big\|_{\Lp}\\
&\le c\Big\|\cM_t\Big(\sum_{\xi\in \cX_j}b_\xi  \ONE_{I_\xi}\Big)(\cdot)\Big\|_{\Lp}
\le c\Big\|\sum_{\xi\in \cX_j}b_\xi \ONE_{I_\xi}(\cdot)\Big\|_{\Lp},
\end{aligned}
\end{equation}
where for the first inequality we used Lemma~\ref{l:weak_inequality}
and for the second Lemma~\ref{lem:disc-max}.
Also, for $\xi\in \cX_j$ let
$\cX_{j+r}(\xi):=\{w\in \cX_{j+r}:I_w\cap I_\xi\ne \emptyset\}$.
Evidently,
$\#\cX_{j+r}(\xi)\le c$.
Then, for $w, \eta\in \cX_{j+r}(\xi)$ we have
$d(w,\eta)\le c(r)2^{-j-r}$ and hence
$$
m_w\le c\frac{m_w}{1+2^{j+r}d(w,\eta)}\le cm^\ast_\eta.
$$
Therefore, for any $\xi\in \cX_j$ and $\eta\in \cX_{j+r}(\xi)$ we have
$b_\xi=\max_{w\in\cX_{j+r}(\xi)}m_w \le cm^\ast_\eta$
and hence
$$
b_\xi\ONE_{I_\xi}\le \sum_{\eta\in \cX_{j+r}(\xi)}  m^\ast_\eta \ONE_{I_\eta}.
$$
Using  this in (\ref{half_shannon}) we get
\begin{equation*}
\begin{aligned}
\Big(\sum_{\xi\in \cX_j}a_\xi^p \mu(I_\xi)\Big)^{1/p}
&\le c\Big\|\sum_{\eta\in \cX_{j+r}}m^\ast_\eta \ONE_{I_\eta}(\cdot)\Big\|_\Lp
\le c\Big\|\cM_t\Big(\sum_{\eta\in \cX_{j+r}}m_\eta\ONE_{I_\eta}\Big)(\cdot)\Big\|_\Lp\\
&\le c\Big\|\sum_{\eta\in \cX_{j+r}}m_\eta  \ONE_{I_\eta}(\cdot)\Big\|_\Lp
\le c\norm{P}_\Lp,
\end{aligned}
\end{equation*}
which completes the proof.
$\qed$

\end{document}